# LIMIT VELOCITY AND ZERO–ONE LAWS FOR DIFFUSIONS IN RANDOM ENVIRONMENT

By Laurent Goergen

*ETH Zurich*

We prove that multidimensional diffusions in random environment have a limiting velocity which takes at most two different values. Further, in the two-dimensional case we show that for any direction, the probability to escape to infinity in this direction equals either zero or one. Combined with our results on the limiting velocity, this implies a strong law of large numbers in two dimensions.

**1. Introduction.** Over the last 25 years, diffusions in a random medium have been the object of many studies. They came as a natural way to generalize homogenization in a periodic medium and model disorder at a microscopic scale; see [3, 16]. In spite of a large literature (see, e.g., [9, 10, 11, 13, 14, 15, 17, 18, 19, 20, 23, 25, 26, 32, 35]), only partial results are known on such basic questions as zero–one laws, recurrence–transience, the law of large numbers and central limit theorems.

The method of the environment viewed from the particle has been a powerful tool in the study of diffusions in a random medium, but many examples fall outside its scope. Recently in the discrete setting, other methods, for instance exploring renewal-type arguments, have contributed to a revival of the subject; see [5, 6, 21, 29, 30, 31, 33, 34, 36, 37, 38]. It is natural, but not straightforward, to try to build up on these ideas and make progress in the continuous framework. This approach has proved successful notably in the ballistic case, that is, when the diffusion has a nonvanishing limiting velocity; see, for instance, [12, 25, 26]. The present article follows a similar endeavor. We prove in the general framework of diffusions in a random environment (see below) the existence of a limiting velocity as well as certain zero–one laws. Corresponding results are known in the discrete framework; see [22, 34, 37, 38]. Our work is closer in spirit to the last two references. It also draws on the renewal structure constructed by Shen [26] which is more intricate than its discrete counterpart in [33].

Before we discuss our results any further, we first describe the model. The random environment is specified by a probability space $(\Omega, \mathcal{A}, \mathbb{P})$ on which acts a jointly measurable group $\{t_x; x \in \mathbb{R}^d\}$ of $\mathbb{P}$-preserving transformations, with $d \geq 1$. The diffusion matrix and the drift of the diffusion in random environment are stationary functions $a(x,\omega)$, $b(x,\omega)$, $x \in \mathbb{R}^d$, $\omega \in \Omega$, with respective values in the space of nonnegative $d \times d$ matrices and in $\mathbb{R}^d$, that is,

$$
\begin{aligned}
a(x+y,\omega) &= a(x, t_y\omega), \\
b(x+y,\omega) &= b(x, t_y\omega) \qquad \text{for } x,y \in \mathbb{R}^d, \omega \in \Omega.
\end{aligned}
\tag{1.1}
$$

We assume that these functions are bounded and uniformly Lipschitz, that is, there is a $\bar{K} > 1$, such that for $x, y \in \mathbb{R}^d, \omega \in \Omega$,

$$
\begin{aligned}
|b(x,\omega)| + |a(x,\omega)| &\leq \bar{K}, \\
|b(x,\omega) - b(y,\omega)| + |a(x,\omega) - a(y,\omega)| &\leq \bar{K}|x-y|,
\end{aligned}
\tag{1.2}
$$

where $|\cdot|$ denotes the Euclidean norm for vectors and matrices. Further we assume that the diffusion matrix is uniformly elliptic, that is, there is a $\nu > 1$ such that for all $x, y \in \mathbb{R}^d, \omega \in \Omega$:

$$
\frac{1}{\nu}|y|^2 \leq y \cdot a(x,\omega)y \leq \nu|y|^2.
\tag{1.3}
$$

The coefficients $a, b$ satisfy a condition of finite range dependence: for $A \subset \mathbb{R}^d$, we define

$$
\mathcal{H}_A = \sigma(a(x,\cdot), b(x,\cdot); x \in A),
\tag{1.4}
$$

and assume that for some $R > 0$,

(1.5)   $\mathcal{H}_A$ and $\mathcal{H}_B$ are independent under $\mathbb{P}$ whenever $d(A,B) \geq R$,

where $d(A,B)$ is the mutual Euclidean distance between $A$ and $B$. With the above regularity assumptions on $a$ and $b$, for any $\omega \in \Omega$, $x \in \mathbb{R}^d$, the martingale problem attached to $x$ and the operator

$$
\mathcal{L}_\omega = \tfrac{1}{2} \sum_{i,j=1}^d a_{ij}(\cdot,\omega)\partial^2_{ij} + \sum_{i=1}^d b_i(\cdot,\omega)\partial_i
\tag{1.6}
$$

is well posed; see [28] or [2], page 130. The corresponding law $P_{x,\omega}$ on $C(\mathbb{R}_+, \mathbb{R}^d)$, unique solution of the above martingale problem, describes the diffusion in the environment $\omega$ and starting from $x$. We write $E_{x,\omega}$ for the expectation under $P_{x,\omega}$ and we denote the canonical process on $C(\mathbb{R}_+, \mathbb{R}^d)$ with $(X_t)_{t\geq 0}$. Observe that $P_{x,\omega}$ is the law of the solution of the stochastic differential equation

(1.7)
$$dX_t = \sigma(X_t, \omega) \, d\beta_t + b(X_t, \omega) \, dt,$$
$$X_0 = x, \qquad P_{x,\omega}\text{-a.s.},$$

where, for instance, $\sigma(\cdot, \omega)$ is the square root of $a(\cdot, \omega)$ and $\beta$ is some $d$-dimensional Brownian motion under $P_{x,\omega}$. The laws $P_{x,\omega}$ are usually called "quenched laws" of the diffusion in random environment. To restore translation invariance, we consider the so-called "annealed laws" $P_x$, $x \in \mathbb{R}^d$, which are defined as semidirect products:

(1.8)
$$P_x \stackrel{\text{def}}{=} \mathbb{P} \times P_{x,\omega}.$$

Of course the Markov property is typically lost under the annealed laws.

The goal of this article is to show the existence of a limiting velocity as well as certain zero–one laws for this process. For any unit vector $l \in \mathbb{R}^d$, denote with

(1.9)
$$A_l = \left\{ \lim_{t \to \infty} l \cdot X_t = +\infty \right\}$$

the event that the diffusion escapes to infinity in direction $l$. We prove a weak zero–one law saying that for any direction $l$, $P_0(A_l \cup A_{-l})$ equals either zero or one; see Proposition 3.6. Then our main result for general dimension $d \geq 1$ (cf. Theorem 3.8) shows the existence of a deterministic unit vector $l_*$ and two deterministic numbers $v_+, v_- \geq 0$, such that

(1.10)
$$\lim_{t \to \infty} \frac{X_t}{t} = (v_+ \mathbb{1}_{A_{l_*}} - v_- \mathbb{1}_{A_{-l_*}}) l_*, \qquad P_0\text{-a.s.}$$

When $d = 2$, we also prove the following stronger zero–one law (cf. Theorem 4.2):

(1.11)
$$\text{for any } l \in S^1, \qquad P_0(A_l) \in \{0, 1\},$$

which together with (1.10) implies the following strong law of large numbers:

(1.12) When $d = 2$, there is a $v \in \mathbb{R}^2$ such that $P_0$-a.s., $\quad \lim_{t \to \infty} \dfrac{X_t}{t} = v.$

In the context of random walks in *ergodic* environments, Zerner and Merkl give in [37] an example, where in the statement corresponding to (1.10) two opposite velocities occur with probability $\frac{1}{2}$ each. This signals that an independence assumption on the environment is of importance for the validity of the zero–one law (1.11) or the law of large numbers (1.12). These questions remain open problems when $d \geq 3$.

To prove (1.10), we consider an arbitrary direction $l$ and proceed differently depending on the value of $P_0(A_l \cup A_{-l})$. In the oscillating case where $P_0(A_l \cup A_{-l}) = 0$, we show in Section 2 that $\lim_{t \to \infty} \frac{l \cdot X_t}{t} = 0$, $P_0$-a.s.; see Corollary 2.6. The argument relies on the fact that for any direction $l \in S^{d-1}$,

(1.13) $\qquad P_0 \left( \limsup_{t \to \infty} \frac{l \cdot X_t}{t} > 0 \right) > 0 \qquad \text{implies } P_0(A_l) > 0;$

see Theorem 2.4. The strategy used to derive (1.13) is similar to the article [38] by Zerner. However, because of finite range dependence and space–time continuity, the arguments are more involved. Nevertheless, we believe that we achieved some simplifications, as our proof avoids infinite products of independent processes (cf. [38] and equation (13) therein). In the context of random walks in a discrete mixing environment, an alternative way to handle the oscillating case can be found in [22].

In order to analyze the case $P_0(A_l \cup A_{-l}) = 1$, we use a renewal structure in the spirit of Shen [26] (see Section 3), and prove that $P_0(A_l) > 0$ implies that on $A_l$, $P_0$-a.s., $\lim_{t \to \infty} \frac{l \cdot X_t}{t} = v_l$. The number $v_l$ is either 0 or expressed in terms of a certain regeneration time $\tau_1$; see (3.41). As in [26], we construct the successive regeneration times $\tau_k$, $k \geq 1$, on an enlarged probability space which is obtained by coupling the diffusion with a suitable sequence of auxiliary i.i.d. Bernoulli variables; see Section 3.1. The quenched measure on the enlarged space, which couples the diffusion to the Bernoulli variables, is denoted with $\hat{P}_{x,\omega}$. In essence, $\tau_1$ is the first time when the trajectory reaches a local maximum in direction $l$, some auxiliary Bernoulli variable takes value 1 and from then on the diffusion never backtracks; see Section 3.2. We generalize the results of Shen to the case where $0 < P_0(A_l) \leq 1$ [instead of assuming $P_0(A_l) = 1$]; see Proposition 3.4 and Theorem 3.5. In the discrete setting, couplings were first used by Zeitouni (cf., e.g., [36], Section 3), with the purpose to overcome the dependence structure of a mixing environment. Another important ingredient for an effective application of the renewal structure is a control on the first moment of $l \cdot X_{\tau_1}$,

(1.14) $\qquad \text{If } P_0(A_l) > 0 \qquad \text{then } \hat{E}_0[l \cdot X_{\tau_1} | D = \infty] < \infty,$

where $\hat{E}_0$ is the expectation under $\mathbb{P} \times \hat{P}_{0,\omega}$ and $\{D = \infty\}$ is the event that the diffusion never backtracks a distance $R$ below its starting point. In the discrete setting, a related result due to Zerner can be found in [36], Lemma 3.2.5. The argument we provide here, however, does not require Blackwell's renewal theorem; see also the comments preceding Proposition 3.7.

In the last section, we prove the zero–one law (1.11) in two dimensions. Our strategy is similar to [37] in the discrete case. We consider two diffusion processes under the law $\mathbb{E}(P_{0,\omega} \times P_{y_L,\omega})$, where $l \cdot y_L \geq 3L$ and $L$ is large. We assume that $P_0(|l \cdot X_t| \to \infty) = 1$ and deduce that the probability of a close encounter of the two diffusions between 0 and $y_L$ vanishes as $L \to \infty$; see Lemma 4.1. This result holds in all dimensions. On the other hand, when $d = 2$, if we assume by contradiction that $P_0(A_l)P_0(A_{-l}) > 0$, we can choose $y_L$ such that for large $L$, the two diffusions intersect "between 0 and $y_L$" with nonvanishing probability; see Theorem 4.2. Then the zero–one law (1.11) follows.

The article is organized as follows. In Section 2 we prove (1.13); see Theorem 2.4. This yields with Corollary 2.6 the main ingredient to prove (1.10) when $P_0(A_l \cup A_{-l}) = 0$. In Section 3 we recall the coupling construction leading to the measures $\hat{P}_{x,\omega}$, define the regeneration times $\tau_k$, $k \geq 1$ (cf. Section 3.2) and develop the theorems describing the renewal structure; see Section 3.3. We also prove a weak zero–one law (cf. Proposition 3.6), as well as (1.14); see Proposition 3.7. Our main result shows for all $d \geq 1$ the existence of a limiting velocity; see (1.10) or Theorem 3.8. In Section 4 we prove the two-dimensional zero–one law (1.11); see Theorem 4.2. In the Appendix we provide for the reader's convenience the proof of a variation of Theorem 2.7 of [26] stated in Lemma 3.3.

*Convention on constants.* Unless otherwise stated, constants only depend on the quantities $d, \bar{K}, \nu, R$. We denote with $c$ positive constants with values changing from place to place and with $c_0, c_1, \ldots$ positive constants with values fixed at their first appearance. Dependence on additional parameters appears in the notation.

**2. Oscillations and null directional speed.** In this section we first introduce some additional notation and then we start with the study of the case, where the trajectory oscillates in some direction $l \in S^{d-1}$. This case corresponds to $P_0[A_l \cup A_{-l}] = 0$ and we will see later that $P_0[A_l \cup A_{-l}]$ equals either zero or one; see Proposition 3.6. The main result is Theorem 2.4: Under the assumption $P_0[\limsup_{t \to \infty} \frac{l \cdot X_t}{t} > 0] > 0$, the trajectories will not backtrack below a certain level with positive probability and with Lemma 2.5, we deduce that $P[A_l] > 0$. It follows then easily that $P_0[A_l \cup A_{-l}] = 0$ implies zero asymptotic speed in the direction $l$ (see Corollary 2.6).

We now introduce some notation used throughout the article. We denote with $\mathbb{N}$ the set of nonnegative integers. The integer part of a real $t \geq 0$ and the smallest integer larger than $t$ are respectively denoted with $\lfloor t \rfloor$ and $\lceil t \rceil$. Let $S^{d-1}$ stand for the Euclidean unit sphere of $\mathbb{R}^d$ and $B(x, r)$ for the open Euclidean ball with radius $r$ centered at $x$. For $a < b$ two reals and $l \in S^{d-1}$, we define

(2.1) $\mathcal{S}_{(a,b)} = \{x \in \mathbb{R}^d; a < x \cdot l < b\}, \qquad \bar{\mathcal{S}}_{(a,b)} = \{x \in \mathbb{R}^d; a \leq x \cdot l \leq b\},$

the open and closed slabs between $a$ and $b$ in the direction $l$. If $A$ is a Borel set of $\mathbb{R}^d$, $|A|$ stands for its Lebesgue measure.

For an open or closed set $A \subset \mathbb{R}^d$, we denote with $H_A = \inf\{t \geq 0; X_t \in A\}$ the entrance time into $A$ and with $T_A = \inf\{t \geq 0; X_t \notin A\}$ the exit time from $A$. We will also use the following stopping times measuring absolute and relative displacements of the trajectory. For $u \in \mathbb{R}$,

$$
\begin{aligned}
T_u &= H_{\{z \in \mathbb{R}^d \,:\, z \cdot l \geq u\}}, \\
\tilde{T}_u &= H_{\{z \in \mathbb{R}^d \,:\, z \cdot l \leq u\}}, \\
T_u^{\text{rel}} &= \inf\{t \geq 0 : l \cdot (X_t - X_0) \geq u\}, \\
\tilde{T}_u^{\text{rel}} &= \inf\{t \geq 0 : l \cdot (X_t - X_0) \leq u\}.
\end{aligned}
\tag{2.2}
$$

We write $(\mathcal{F}_t)_{t \geq 0}$ and $(\theta_t)_{t \geq 0}$ for the canonical right-continuous filtration and for the canonical time-shift on $C(\mathbb{R}_+, \mathbb{R}^d)$, respectively.

We turn now to the construction of the objects appearing in Proposition 2.1. We consider some number $L = 3L' > 3R$ and define the successive times of entrance in $\bar{\mathcal{S}}_{(mL+L', mL+2L')}$ and departure from $\mathcal{S}_{(mL, (m+1)L)}$ (cf. Figure 1): for $m \in \mathbb{N}$,

$$
R_1^{(m)} = H_{\bar{\mathcal{S}}_{(mL+L', mL+2L')}}, \qquad S_1^{(m)} = T_{\mathcal{S}_{(mL, (m+1)L)}} \circ \theta_{R_1^{(m)}} + R_1^{(m)}, \tag{2.3}
$$

and by induction for $k \geq 2$,

$$
R_k^{(m)} = R_1^{(m)} \circ \theta_{S_{k-1}^{(m)}} + S_{k-1}^{(m)}, \qquad S_k^{(m)} = S_1^{(m)} \circ \theta_{S_{k-1}^{(m)}} + S_{k-1}^{(m)}.
$$

We define, for integer $\alpha \geq 2$ (this integer will typically be large in the sequel)

$$
N_\alpha^{(m)} = \sum_{k \geq 1}^{\infty} \mathbb{1}_{\{R_k^{(m)} + 1 \leq S_k^{(m)} < T_{(m+\alpha)L} < \infty\}}, \tag{2.4}
$$

the number of entrances in $\bar{\mathcal{S}}_{(mL+L', mL+2L')}$ after which the trajectory stays at least one time unit in $\mathcal{S}_{(mL, (m+1)L)}$. Moreover, we consider

$$
k_\alpha^{(m)} = \begin{cases} \max\{k \geq 1 : R_k^{(m)} + 1 \leq S_k^{(m)} < T_{(m+\alpha)L}\}, \\ \qquad \text{if } T_{(m+\alpha)L} < \infty \text{ and } \{\cdots\} \neq \varnothing, \\ 0, \qquad \text{otherwise}; \end{cases} \tag{2.5}
$$

$$
h_\alpha^{(m)} = \begin{cases} S_{k_\alpha^{(m)}}^{(m)} - T_{mL}, \\ \qquad \text{if } T_{(m+\alpha)L} < \infty, \text{ with the convention } S_0^{(m)} = T_{mL}, \\ \infty, \qquad \text{otherwise.} \end{cases} \tag{2.6}
$$

The quantity $h_\alpha^{(m)}$ is the time duration, beginning at $T_{mL}$, after which the trajectory does not make "long visits" to the slab $\bar{\mathcal{S}}_{(mL+L', mL+2L')}$ anymore. Note that $h_\alpha^{(m)}$ is nondecreasing in $\alpha$.

Let us give an outline of the steps leading to the main result of this section, that is, Theorem 2.4. In Proposition 2.1 we show that a continuous path $w$ satisfying $\limsup_{t \to \infty} \frac{l \cdot w(t)}{t} > 0$ has the property that there is a large asymptotic fraction of slabs among the $\mathcal{S}_{(mL,(m+1)L)}$, $m \geq 1$, around which the oscillations of $w$ that occur before reaching a level at a distance $\alpha L$ in direction $l$, last only some finite time $h$ independent of $\alpha$. An analogous result for a discrete path is stated in [38], Lemma 3. In the next step, we deduce the existence of an $h > 0$ such that with positive probability the following events, later called $C_m$ [cf. (2.20)], happen with a large asymptotic frequency: on $C_m$, the particle at time $H_{\mathcal{S}_{(mL,(m+1)L)}} + h$ is located in a narrow slab "to the right of" $\mathcal{S}_{(mL,(m+1)L)}$ and then moves to a level at a distance $\alpha L$ without backtracking; see Lemma 2.3. Then we extract the crucial information about the absence of backtracking. In essence for this purpose, we condition each event $C_m$ on the information prior to $H_{\mathcal{S}_{(mL,(m+1)L)}} + h$, and transfer our control on the asymptotic frequency of the $C_m$'s, to a control on the asymptotic mean of the conditional probabilities. This is done with the help of certain martingales and Azuma's inequality; see (2.33). Finally we dominate these conditional probabilities by a sequence of i.i.d. variables under $\mathbb{P}$, apply the law of large numbers and conclude that the probability to never backtrack is positive by letting $\alpha$ tend to infinity. This method bypasses the technique of infinite products of probability spaces in [38] [cf. (13) therein], which is hard to implement in the continuous setting.

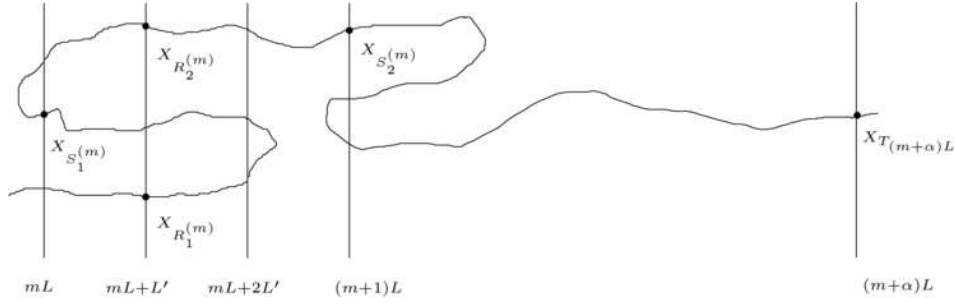

Fig. 1.

PROPOSITION 2.1. *Let $w(\cdot)$ be a continuous path in $\mathbb{R}^d$ starting at $0$ and satisfying $\limsup_{t \to \infty} \frac{l \cdot w(t)}{t} > 0$; then there exists an integer $h \geq 1$ such that for all integers $\alpha \geq 2$,*

$$(2.7) \qquad \limsup_{M \to \infty} \frac{1}{M+1} \sum_{m=0}^{M} \mathbb{1}_{\{h_\alpha^{(m)}(w) \leq h\}} \geq \frac{1}{3}.$$

PROOF. We choose $\delta > 0$ such that $\limsup_{t \to \infty} \frac{l \cdot w(t)}{t} \geq \delta$. There is a sequence $(t_k)_{k \geq 1}$ in $\mathbb{R}_+$ tending to infinity such that $l \cdot w(t_k) > \delta t_k$. Thus, for all $\alpha \geq 2$,

(2.8) $T_{(\delta/2)t_k + \alpha L}(w) \leq t_k$ for all large enough $k$ (depending on $\alpha$).

(For the sake of simplicity, we will drop $w$ from the notation.) If we choose $M_k$ integer such that $M_k L \leq \frac{\delta}{2} t_k \leq (M_k + 1)L, k \geq 1$, (2.8) implies that for all integers $\alpha \geq 2$,

$$T_{(M_k+\alpha)L} \leq \frac{2(M_k+1)L}{\delta} \qquad \text{for all large enough } k \text{ (depending on } \alpha\text{)}. \tag{2.9}$$

If $R_k^{(m)} + 1 \leq S_k^{(m)} < T_{(m+\alpha)L}$ and since $T_{(m+\alpha)L}$ is finite for all $m$, the path $w$ spends at least one unit of time entirely in the slab $\mathcal{S}_{(mL,(m+1)L)}$ before reaching level $(m+\alpha)L$. Hence, for all $k$ large enough, we deduce from (2.9) that [cf. (2.4) for the notation]

$$\sum_{m=0}^{M_k} N_\alpha^{(m)} \leq T_{(M_k+\alpha)L} \leq \frac{2(M_k+1)L}{\delta} \tag{2.10}$$

and

$$\sum_{m=0}^{M_k} h_\alpha^{(m)} \leq \sum_{j=0}^{\alpha-1} \sum_{\substack{m \bmod \alpha = j \\ m \leq M_k}} (T_{(m+\alpha)L} - T_{mL})$$

$$\leq \alpha T_{(M_k+\alpha)L} \leq \frac{2\alpha(M_k+1)L}{\delta}, \tag{2.11}$$

for all large enough $k$. Assume now that (2.7) with real $h$ does not hold, that is: for all $h \geq 1$, there is an integer $\alpha \geq 2$ such that $\limsup_{M \to \infty} \frac{1}{M+1} \sum_{m=0}^{M} \mathbb{1}_{\{h_\alpha^{(m)} \leq h\}} < \frac{1}{3}$. We can construct inductively $h_0 = 1$, $\alpha_1 \geq 2, h_i = \frac{6\alpha_i L}{\delta}, \alpha_{i+1} > \alpha_i$ using that $h_\alpha^{(m)}$ is nondecreasing in $\alpha$, such that

$$\text{for all } i \geq 1 \qquad \limsup_{M \to \infty} \frac{1}{M+1} \sum_{m=0}^{M} \mathbb{1}_{\{h_{\alpha_{i+1}}^{(m)} \leq h_i\}} < \frac{1}{3}. \tag{2.12}$$

On the other hand, (2.11) and the choice $h_i = \frac{6\alpha_i L}{\delta}$ imply that

$$\text{for all } i \geq 1 \qquad \limsup_{k \to \infty} \frac{1}{M_k+1} \sum_{m=0}^{M_k} \mathbb{1}_{\{h_{\alpha_i}^{(m)} > h_i\}} \leq \frac{1}{3}. \tag{2.13}$$

Observe that for all $i, k \geq 1$,

$$1 \leq \frac{1}{M_k+1} \sum_{m=0}^{M_k} \mathbb{1}_{\{h_{\alpha_{i+1}}^{(m)} \leq h_i\}} + \mathbb{1}_{\{h_{\alpha_i}^{(m)} > h_i\}} + \mathbb{1}_{\{h_{\alpha_i}^{(m)} < h_{\alpha_{i+1}}^{(m)}\}}.$$

This inequality together with (2.12) and (2.13) yields

$$\text{for all } i \geq 1 \qquad \frac{1}{3} \leq \frac{1}{M_k+1} \sum_{m=0}^{M_k} \mathbb{1}_{\{h_{\alpha_i}^{(m)} < h_{\alpha_{i+1}}^{(m)}\}} \tag{2.14}$$

for all large enough $k$.

If $h_{\alpha_i}^{(m)} < h_{\alpha_{i+1}}^{(m)}$, the trajectory, after reaching level $(m+\alpha_i)L$, has to return to the slab $\bar{\mathcal{S}}_{(mL+L',mL+2L')}$ and stay in the slab $\mathcal{S}_{(mL,(m+1)L)}$ for at least one unit of time, all this before reaching level $(m+\alpha_{i+1})L$. Therefore we see that

$$\mathbb{1}_{\{h_{\alpha_i}^{(m)} < h_{\alpha_{i+1}}^{(m)}\}} \leq N_{\alpha_{i+1}}^{(m)} - N_{\alpha_i}^{(m)},$$

and hence for arbitrary $i_0 \geq 1$ and large $k$, we obtain

$$\frac{i_0}{3} \stackrel{(2.14)}{\leq} \frac{1}{M_k+1} \sum_{m=0}^{M_k} \sum_{i=1}^{i_0} \mathbb{1}_{\{h_{\alpha_i}^{(m)} < h_{\alpha_{i+1}}^{(m)}\}} \leq \frac{1}{M_k+1} \sum_{m=0}^{M_k} N_{\alpha_{i_0}+1}^{(m)} \stackrel{(2.10)}{\leq} \frac{2L}{\delta}, \quad (2.14)$$

a contradiction. We have thus proved the existence of a real $h \geq 1$ such that (2.7) holds. By monotonicity, we can increase $h$ to be an integer. $\square$

The next lemma comes as a preparation for the main result of this section, namely Theorem 2.4. If $S$ is any stopping time, we write $S_k$, $k \geq 0$, for the iterates of $S$, namely,

$$(2.15) \qquad S_0 = 0, \qquad S_1 = S, \qquad S_{k+1} = S \circ \theta_{S_k} + S_k \leq \infty.$$

LEMMA 2.2. *For every $k \geq 1$, let $U^k$ be an $(\mathcal{F}_t)_{t \geq 0}$-stopping time and $\Delta_k \in \mathcal{F}_{U^k}$. Denote with $U_m^k$, $m \geq 0$, the iterates of $U^k$. If there exist numbers $\gamma_1, \gamma_2 > 0$, such that*

$$(2.16) \qquad \text{for all } k \geq 1, x \in \mathbb{R}^d, \omega \in \Omega : P_{x,\omega}(\Delta_k) \leq \gamma_1 e^{-\gamma_2 k},$$

*then for each $\varepsilon > 0$, there is a $k_0(\varepsilon, \gamma_1, \gamma_2) \geq 1$, such that*

$$(2.17) \quad \text{for } k \geq k_0, P_0\text{-a.s.,} \qquad \limsup_{M \to \infty} \frac{1}{M+1} \sum_{m=0}^{M} \mathbb{1}_{\Delta_k} \circ \theta_{U_m^k} \leq \varepsilon,$$

*with the convention that $\mathbb{1}_{\Delta_k} \circ \theta_{U_m^k} = 0$, if $U_m^k = \infty$.*

PROOF. Note that $\mathbb{1}_{\Delta_k} \circ \theta_{U_m^k}$ is $\mathcal{F}_{U_{m+1}^k}$-measurable. The strong Markov property yields for $M \geq 1, k \geq 1$

$$E_{x,\omega}\left[\exp\left(\sum_{m=0}^{M} \mathbb{1}_{\Delta_k} \circ \theta_{U_m^k}\right), U_M^k < \infty\right]$$

$$(2.18) \quad = E_{x,\omega}\left[\exp\left(\sum_{m=0}^{M-1} \mathbb{1}_{\Delta_k} \circ \theta_{U_m^k}\right) E_{X_{U_M^k},\omega}[\exp(\mathbb{1}_{\Delta_k})], U_M^k < \infty\right]$$

$$= E_{x,\omega}\left[\exp\left(\sum_{m=0}^{M-1} \mathbb{1}_{\Delta_k} \circ \theta_{U_m^k}\right) ((e-1) P_{X_{U_M^k},\omega}[\Delta_k] + 1), U_M^k < \infty\right].$$

Using (2.16) and iteration, we obtain for $M \geq 1$, $k \geq 1$,

$$E_{x,\omega}\left[\exp\left(\sum_{m=0}^{M} \mathbb{1}_{\Delta_k} \circ \theta_{U_m^k}\right), U_M^k < \infty\right] \leq (\gamma_1 e^{-\gamma_2 k}(e-1) + 1)^{M+1}.$$

Therefore, using Chebyshev's inequality, we find

$$P_{x,\omega}\left[\frac{1}{M}\sum_{m=0}^{M-1} \mathbb{1}_{\Delta_k} \circ \theta_{U_m^k} > \varepsilon\right]$$
$$\leq e^{-\varepsilon M}(\gamma_1 e^{-\gamma_2 k}(e-1) + 1)^M \leq e^{M(-\varepsilon + (e-1)\gamma_1 e^{-\gamma_2 k})}.$$

If $k$ is large enough, the argument of the exponential becomes negative and our claim follows from Borel–Cantelli's lemma. □

In the next lemma with two successive reduction steps, we replace $\{h_\alpha^{(m)} \leq h\}$ appearing in (2.7) by an event $C_m$ that has the following meaning: At the stopping time $T_{mL} + h_0$, for some $h_0 \geq 1$, the position of the diffusion is located in the slab $\mathcal{S}_{(mL+2L',(m+K)L)}$ and after this stopping time, the trajectory reaches level $(m+\alpha)L$ without going below level $mL + 2L'$; see (2.20).

LEMMA 2.3. *Assume that $P_0(\limsup_{t\to\infty} \frac{l\cdot X_t}{t} > 0) > 0$; then there exist an integer $h_0 \geq 1$ and constants $L = 3L' > 3R, K = K(h_0) \geq 1$, such that*

$$(2.19) \qquad P_0\left[\inf_{\alpha \geq K} \limsup_{M\to\infty} \frac{1}{M+1} \sum_{m=0}^{M} \mathbb{1}_{C_m} \geq \frac{1}{12}\right] > 0,$$

*where*

$$(2.20) \quad \begin{aligned} C_m &= \{X_{T_{mL}+h_0} \in \mathcal{S}_{(mL+2L',(m+K)L)}\} \cap \{T_{mL} + h_0 < T_{(m+\alpha)L}\} \\ &\quad \cap \theta_{T_{mL}+h_0}^{-1}\{\tilde{T}_{mL+2L'} > T_{(m+\alpha)L}\} \qquad \text{for } m \geq 0. \end{aligned}$$

PROOF. With our assumption, Proposition 2.1 yields that for some integer $h_0$

$$(2.21) \qquad P_0\left[\inf_{\alpha \geq 2} \limsup_{M\to\infty} \frac{1}{M+1} \sum_{m=0}^{M} \mathbb{1}_{\{h_\alpha^{(m)} \leq h_0\}} \geq \frac{1}{3}\right] > 0.$$

In a first reduction step, we want to keep only those slabs $\mathcal{S}_{(mL,(m+1)L)}$, where after time $T_{mL} + h_0$ and before reaching level $(m+\alpha)L$, the paths do not return to the inner part $\bar{\mathcal{S}}_{(mL+L',mL+2L')}$ of the slab. More precisely we claim that if $L' > R$ is large enough, then we obtain from (2.21)

$$(2.22)\ P_0\left[\inf_{\alpha \geq 2} \limsup_{M\to\infty} \frac{1}{M+1} \sum_{m=0}^{M} \mathbb{1}_{\{h_\alpha^{(m)} \leq h_0, R' \circ \theta_{T_{mL}} + T_{mL} > T_{(m+\alpha)L}\}} \geq \frac{1}{6}\right] > 0,$$

where we used the notation

(2.23) $$R' = \inf\{t \geq h_0 : l \cdot (X_t - X_0) \in [L', 2L']\}.$$

Indeed, consider for fixed $\alpha \geq 2$, $m \geq 0$, a trajectory $w$ starting in 0 and satisfying $h_\alpha^{(m)} \leq h_0$. If $w$ visits the inner slab $\bar{\mathcal{S}}_{(mL+L',mL+2L')}$ between time $T_{mL} + h_0$ and $T_{(m+\alpha)L}$, then it must exist from the outer slab $\mathcal{S}_{(mL,(m+1)L)}$ within time 1 as otherwise $h_\alpha^{(m)}$ becomes larger than $h_0$. Note also that by definition, $h_\alpha^{(m)} \leq h_0$ implies $T_{(m+\alpha)L} < \infty, P_0$-a.s. Hence $P_0$-a.s. we have

$$\{h_\alpha^{(m)} \leq h_0, R' \circ \theta_{T_{mL}} + T_{mL} \leq T_{(m+\alpha)L}\} \subset \theta_{T_{mL}}^{-1}(\Delta_{L',\alpha}) \cap \{T_{mL} < \infty\},$$

where we defined $\Delta_{L',\alpha} = \{\sup_{s \leq 1} |X_s - X_0| \circ \theta_{R'} \geq L', R' \leq T_{\alpha L}^{\text{rel}} < \infty\}$. By the Markov property and Bernstein's inequality (see [2], Proposition 8.1, page 23), we have for all $x \in \mathbb{R}^d$, $\omega \in \Omega$,

(2.24) $$P_{x,\omega}[\Delta_{L',\alpha}] \leq E_{x,\omega}\left[R' < \infty, P_{X_{R'},\omega}\left[\sup_{s \leq 1}|X_s - X_0| \geq L'\right]\right]$$
$$\leq c_1 e^{-c_2 L'^2}.$$

We decompose the indicator-function $\mathbb{1}_{\{h_\alpha^{(m)} \leq h_0\}}$ appearing in (2.21) as follows:

(2.25) $$\limsup_{M \to \infty} \frac{1}{M+1} \sum_{m=0}^{M} \mathbb{1}_{\{h_\alpha^{(m)} \leq h_0\}}$$
$$\leq \limsup_{M \to \infty} \frac{1}{M+1} \sum_{m=0}^{M} \mathbb{1}_{\{h_\alpha^{(m)} \leq h_0\}} \mathbb{1}_{\{R' \circ \theta_{T_{mL}} + T_{mL} > T_{(m+\alpha)L}\}}$$
$$+ \limsup_{M \to \infty} \frac{1}{M+1} \sum_{m=0}^{M} \mathbb{1}_{\Delta_{L',\alpha}} \circ \theta_{T_{mL}}.$$

In order to apply Lemma 2.2 to the last term of (2.25), since $\theta_{T_{mL}}^{-1} \Delta_{L',\alpha} \in \mathcal{F}_{T_{(\alpha+m)L}}$, for $m \geq 0$, we rewrite the sum in the last term as a double sum running over all residue classes modulo $\alpha$ and obtain as an upper bound

(2.26) $$\frac{1}{\alpha} \sum_{j=0}^{\alpha-1} \limsup_{M \to \infty} \frac{1}{M+1} \sum_{m=0}^{M} \mathbb{1}_{\Delta_{L',\alpha}} \circ \theta_{T_{(m\alpha+j)3L'}}.$$

With (2.24), we can apply Lemma 2.2 for every $j = 0, \ldots, \alpha - 1$. The parameter $L'$ is chosen integer and plays the role of $k$ in the lemma. Moreover we, respectively, substitute $T_{\alpha 3L'}^{\text{rel}} \circ \theta_{T_{jL}} + T_{jL}$ and $\theta_{T_{j3L'}}^{-1}(\Delta_{L',\alpha}) \cap \{T_{j3L'} < \infty\}$ for $U^k$ and $\Delta_k$, and use $\varepsilon = \frac{1}{6}$. Note that $P_0$-a.s., for $m \geq 1$, the $m$th iterate of $T_{\alpha 3L'}^{\text{rel}} \circ \theta_{T_{jL}} + T_{jL}$ is $T_{(m\alpha+j)3L'}$. As the lower bound for $L'$ provided by Lemma 2.2 only depends on the constants $c_1, c_2$ in (2.24), there exists a constant $L' > R$, such that for all $\alpha \geq 2$ and all $j = 0, \ldots, \alpha - 1$, we obtain

$$\limsup_{M \to \infty} \frac{1}{M+1} \sum_{m=0}^{M} \mathbb{1}_{\Delta_{L',\alpha}} \circ \theta_{T_{(m\alpha+j)L}} < \frac{1}{6}, \qquad P_0\text{-a.s.}$$

Hence the last term of (2.25) is $P_0$-a.s. for all integers $\alpha \geq 2$ smaller than $\frac{1}{6}$. In view of (2.21), this estimate proves the claim (2.22) of the first reduction step.

In the second reduction step, we would like to keep only those slabs where the trajectory stays in a big ball during time $h_0$ after $T_{mL}$. We claim that there exists a constant $K = K(h_0, L) \geq 1$, such that

$$(2.27) \quad 0 < P_0 \Bigg[ \inf_{\alpha > K} \limsup_{M \to \infty} \frac{1}{M+1} \sum_{m=0}^{M} \mathbb{1}_{\{h_\alpha^{(m)} \leq h_0, R' \circ \theta_{T_{mL}} + T_{mL} > T_{(m+\alpha)L}\}}$$

$$\times \mathbb{1}_{\{\sup_{s \leq h_0} |X_s - X_0| \circ \theta_{T_{mL}} < KL\}} \geq \frac{1}{12} \Bigg].$$

Indeed, define $\Delta'_k = \{\sup_{s \leq h_0} |X_s - X_0| \geq kL\}$, for $k \geq 1$. From Bernstein's inequality (see [2], Proposition 8.1, page 23), there exist positive constants $c_3(h_0), c_4(h_0)$, such that

$$(2.28) \quad P_{x,\omega}[\Delta'_k] \leq c_3 e^{-c_4 (kL)^2} \qquad \text{for all } x \in \mathbb{R}^d, \omega \in \Omega.$$

As before we decompose the indicator-function appearing in (2.22) according to $\theta_{T_{mL}}^{-1}(\Delta'_k{}^c)$ and $\theta_{T_{mL}}^{-1}(\Delta'_k)$. In view of an application of Lemma 2.2, since $\theta_{T_{mL}}^{-1} \Delta'_k \in \mathcal{F}_{T_{(k+m)L}}$, $m \geq 0$, we rewrite the sum $\sum_{m=0}^{M} \mathbb{1}_{\Delta'_k} \circ \theta_{T_{mL}}$ as a double sum running over all the residue classes modulo $k$. As a result

$$(2.29) \quad \begin{aligned} &\limsup_{M \to \infty} \frac{1}{M+1} \sum_{m=0}^{M} \mathbb{1}_{\Delta'_k} \circ \theta_{T_{mL}} \\ &\leq \frac{1}{k} \sum_{j=0}^{k-1} \limsup_{M \to \infty} \frac{1}{M+1} \sum_{m=0}^{M} \mathbb{1}_{\Delta'_k} \circ \theta_{T_{(mk+j)L}}. \end{aligned}$$

With (2.28), we apply Lemma 2.2 for every $j = 0, \ldots, k-1$: $\theta_{T_{jL}}^{-1}(\Delta'_k) \cap \{T_{jL} < \infty\}$ and $T_{kL}^{\text{rel}} \circ \theta_{T_{jL}} + T_{jL}$ play the role of $\Delta_k$ and $U^k$ and we choose $\varepsilon = \frac{1}{12}$ in Lemma 2.2. $P_0$-a.s., for $m \geq 1$, the $m$th iterate of $T_{kL}^{\text{rel}} \circ \theta_{T_{jL}} + T_{jL}$ is $T_{(mk+j)L}$. Hence there is a constant $K = K(h_0, L) \geq 1$, such that for $k = K$ the left-hand side of (2.29) is $P_0$-a.s. smaller than $\frac{1}{12}$. This proves (2.27). We conclude the proof by noting that $\{h_\alpha^{(m)} \leq h_0, \sup_{s \leq h_0} |X_s - X_0| \circ \theta_{T_{mL}} < KL, R' \circ \theta_{T_{mL}} > T_{(m+\alpha)L} - T_{mL}\}$ is $P_0$-a.s. included in $C_m$. $\square$

THEOREM 2.4. $(d \geq 1)$. *For any $l \in S^{d-1}$, $P_0(\limsup_{t \to \infty} \frac{l \cdot X_t}{t} > 0) > 0$ implies that there exists a number $r_0 > 0$ such that $P_0(\tilde{T}^{\mathrm{rel}}_{-r_0} = \infty) > 0$ and as a consequence $P_0(A_l) > 0$.*

(Note that $r_0$ only depends on the quantities $h_0, L, K$ from Lemma 2.3.)

PROOF. Let us briefly outline the argument: We would like to apply the law of large numbers to the sum $\sum_{m=0}^{M} \mathbb{1}_{C_m}$ appearing in (2.19), but the dependence structure of the sequence $(C_m)_{m \geq 0}$ seems to be complicated. Therefore we will replace this sequence by one that is i.i.d. with respect to $\mathbb{P}$. This will be achieved by constructing an appropriate martingale and using Azuma's inequality.

We pick $L = 3L', K$ and $h_0$ as in Lemma 2.3. We introduce the following filtrations: for integer $\alpha > K$ and $j = 0, \ldots, \alpha - 1$,

$$\begin{aligned}
\mathcal{G}^j_m &= \mathcal{F}_{T_{(m\alpha+j)L}}, & m \geq 1, \quad \mathcal{G}^j_0 = \mathcal{F}_0; \\
\tilde{\mathcal{G}}^j_m &= \mathcal{F}_{T_{((m+1)\alpha+j)L} \wedge (T_{(m\alpha+j)L} + h_0)}, & m \geq 0.
\end{aligned} \tag{2.30}$$

Recall the definition of $C_m$ (2.20) and observe that for $\alpha > K$, $j = 0, \ldots, \alpha - 1$, $m \geq 0$:

$$C_{m\alpha+j} \in \mathcal{G}^j_{m+1} \quad \text{and} \quad \mathcal{G}^j_m \subset \tilde{\mathcal{G}}^j_m \subset \mathcal{G}^j_{m+1}, \tag{2.31}$$

because $T_{(m\alpha+j)L} \leq T_{((m+1)\alpha+j)L} \wedge (T_{(m\alpha+j)L} + h_0) \leq T_{((m+1)\alpha+j)L}$.

We define for $j = 0, \ldots, \alpha - 1$ and $n \geq 1$

$$M^j_n = \sum_{m=0}^{n-1} \mathbb{1}_{C_{m\alpha+j}} - E_{0,\omega}[\mathbb{1}_{C_{m\alpha+j}} | \tilde{\mathcal{G}}^j_m], \qquad M^j_0 = 0. \tag{2.32}$$

By (2.31), $M^j_n$ is $\mathcal{G}^j_n$-measurable, for $n \geq 0$, and integrable. It is a $\mathcal{G}^j_n$ martingale under $P_{0,\omega}$, for any $\omega \in \Omega$, because $P_{0,\omega}$-a.s., for $n \geq 1$ we have

$$E_{0,\omega}[\mathbb{1}_{C_{(n-1)\alpha+j}} | \mathcal{G}^j_{n-1}] = E_{0,\omega}[E_{0,\omega}[\mathbb{1}_{C_{(n-1)\alpha+j}} | \tilde{\mathcal{G}}^j_{n-1}] | \mathcal{G}^j_{n-1}].$$

Since $M^j_m$ has bounded increments Azuma's inequality (see [1], Theorem 2.1) applies and we find

$$P_{0,\omega}\left[\lim_{n \to \infty} \frac{M^j_n}{n} = 0\right] = 1 \qquad \text{for all } j = 0, \ldots, \alpha - 1.$$

Hence for any $\omega \in \Omega$, $P_{0,\omega}$-a.s., for all $\alpha > K$ and $j = 0, \ldots, \alpha - 1$,

$$\limsup_{M \to \infty} \frac{1}{M+1} \sum_{m=0}^{M} \mathbb{1}_{C_{m\alpha+j}} = \limsup_{M \to \infty} \frac{1}{M+1} \sum_{m=0}^{M} E_{0,\omega}[\mathbb{1}_{C_{m\alpha+j}} | \tilde{\mathcal{G}}^j_m]. \tag{2.33}$$

The strong Markov property yields that $P_{0,\omega}$-a.s., for $\alpha > K$, $j = 0, \ldots, \alpha - 1$ and $m \geq 0$

$$E_{0,\omega}[\mathbb{1}_{C_{m\alpha+j}} | \tilde{\mathcal{G}}_m^j] = \mathbb{1}_{\{X_{T_{(m\alpha+j)L}+h_0} \in I_{m,j}, T_{(m\alpha+j)L}+h_0 < T_{((m+1)\alpha+j)L}\}}$$

(2.34)
$$\times P_{X_{T_{(m\alpha+j)L}+h_0},\omega}[\tilde{T}_{(m\alpha+j)L+2L'} > T_{((m+1)\alpha+j)L}]$$

$$\leq \sup_{y \in I_{m,j}} P_{y,\omega}[\tilde{T}_{(m\alpha+j)L+2L'} > T_{((m+1)\alpha+j)L}],$$

where $I_{m,j} \stackrel{\text{def}}{=} \bar{\mathcal{S}}_{((m\alpha+j)L+2L',(m\alpha+j+K)L)}$.

Therefore, (2.19) together with (2.33) and (2.34) imply that

(2.35)
$$\mathbb{P}\left[\inf_{\alpha > K} \frac{1}{\alpha} \sum_{j=0}^{\alpha-1} \limsup_{M \to \infty} \frac{1}{M+1} \right.$$
$$\left. \times \sum_{m=0}^{M} \sup_{y \in I_{m,j}} P_{y,\omega}[\tilde{T}_{(m\alpha+j)L+2L'} > T_{((m+1)\alpha+j)L}] > \frac{1}{12}\right] > 0.$$

With respect to $\mathbb{P}$, the variables $f_{m,j} \stackrel{\text{def}}{=} \sup_{y \in I_{m,j}} P_{y,\omega}[\tilde{T}_{(m\alpha+j)L+2L'} > T_{((m+1)\alpha+j)L}]$, $m \geq 0$, are i.i.d. for every $j = 0, \ldots, \alpha - 1$. Indeed the respective slabs $\bar{\mathcal{S}}_{((m\alpha+j)L+2L',((m+1)\alpha+j)L)}$ as $m$ varies are separated by at least $2L' > R$, and one applies (1.5), as well as translation invariance. Hence, from the law of large numbers and from (2.35) we deduce that

$$\inf_{\alpha > K} \mathbb{E}\left[\sup_{y \in I_{0,0}} P_{y,\omega}[\tilde{T}_{2L'} > T_{\alpha L}]\right] > \tfrac{1}{12},$$

and by dominated convergence for $\alpha \to \infty$

(2.36) $\quad \tfrac{1}{12} < \mathbb{E}\left[\sup_{y \in I_{0,0}} P_{y,\omega}[\tilde{T}_{2L'} = \infty]\right] \leq \mathbb{E}\left[\sup_{y \in I_{0,0}} P_{y,\omega}[\tilde{T}^{\text{rel}}_{-(KL-2L')} = \infty]\right].$

We define $\frac{r_0}{2} = KL - 2L'$. Assume now that $P_0[\tilde{T}^{\text{rel}}_{-r_0} = \infty] = 0$. It follows from Fubini's theorem that there is a $\mathbb{P}$-null set $\Gamma \subset \Omega$, such that

(2.37)
$$\text{for } \omega \in \Gamma^c \quad P_{x,\omega}[\tilde{T}^{\text{rel}}_{-r_0} = \infty] = 0,$$
$$\text{except on a Lebesgue-negligible subset of } \mathbb{R}^d.$$

But then, for any $y \in \mathbb{R}^d, \omega \in \Omega$,

(2.38)
$$P_{y,\omega}[\tilde{T}^{\text{rel}}_{-r_0/2} = \infty]$$
$$= P_{y,\omega}\left[\tilde{T}^{\text{rel}}_{-r_0/2} = \infty, \sup_{s \leq 1/n} |X_s - X_0| \leq \frac{r_0}{4}\right]$$

$$+ P_{y,\omega}\left[\tilde{T}^{\text{rel}}_{-r_0/2} = \infty, \sup_{s \leq 1/n} |X_s - X_0| > \frac{r_0}{4}\right]$$

$$\leq P_{y,\omega}[\tilde{T}^{\text{rel}}_{-r_0} \circ \theta_{1/n} = \infty] + P_{y,\omega}\left[\sup_{s \leq 1/n} |X_s - X_0| > \frac{r_0}{4}\right].$$

By the Markov property, the first term on the right-hand side equals $\int_{\mathbb{R}^d} p^\omega(y, x, \frac{1}{n}) P_{x,\omega}[\tilde{T}^{\text{rel}}_{-r_0} = \infty] \, dx$, where $p^\omega(y, \cdot, \frac{1}{n})$ denotes the transition density of the diffusion starting in $y$ in the environment $\omega$ at time $1/n$ with respect to the Lebesgue measure. This density exists under the assumptions (1.2), (1.3); see [4], Theorem 4.5. Hence using (2.37), this term equals 0 for all $\omega \in \Gamma^c$, $y \in \mathbb{R}^d$. The second term on the right-hand side of (2.38) converges to 0 as $n \to \infty$ by continuity of the trajectories. And so it would follow that for all $\omega \in \Gamma^c, y \in \mathbb{R}^d$, $P_{y,\omega}[\tilde{T}^{\text{rel}}_{-r_0/2} = \infty] = 0$. But this contradicts (2.36) and hence

$$(2.39) \qquad P_0[\tilde{T}^{\text{rel}}_{-r_0} = \infty] > 0.$$

To show that $P_0(A_l) > 0$, we need the following useful lemma:

LEMMA 2.5. *Consider $l$ in $S^{d-1}$. For $u, v \in \mathbb{R}, u < v$, define the stopping times $\beta^{l,u} = \inf\{t \geq 1; l \cdot X_t \geq u\}$ and denote their iterates with $\beta^{l,u}_k$, $k \geq 0$. Then one has for all $x \in \mathbb{R}^d$, $\omega \in \Omega$*

$$(2.40) \qquad P_{x,\omega}[\beta^{l,u}_k < \infty, \text{ for all } k \geq 0 \text{ and } T_v = \infty] = 0.$$

PROOF. By the support theorem (see [2], page 25), there is a constant $c = c(v - u) > 0$ such that for all $x \in \bar{\mathcal{S}}_{(u,v)}$, $\omega \in \Omega : P_{x,\omega}[T_v \leq 1] > c$. Then the Markov property shows that for all $x \in \mathbb{R}^d$, $\omega \in \Omega$, $k \geq 1$,

$$P_{x,\omega}[\beta^{l,u}_k < \infty, T_v = \infty] \leq P_{x,\omega}[0 \leq \beta^{l,u}_k < T_v]$$
$$\leq P_{x,\omega}[0 \leq \beta^{l,u}_{k-1} < T_v, T_v \circ \theta_{\beta^{l,u}_{k-1}} > 1]$$
$$\leq (1-c) P_{x,\omega}[0 \leq \beta^{l,u}_{k-1} < T_v].$$

After iteration and letting $k$ tend to infinity, we obtain the claim. □

We are now ready to finish the proof of Theorem 2.4. We observe that for any $v > 0$

$$(2.41) \qquad \{\tilde{T}_{-v} = \infty\} \subset A_l, \qquad P_0\text{-a.s.}$$

Indeed, we have in view of Lemma 2.5 with $-l$ in the role of $l$

$$P_0[A^c_l, \tilde{T}_{-v} = \infty] = P_0[\text{for some } u \in \mathbb{Z}, u < v : \beta^{-l,u}_k < \infty \text{ and } \tilde{T}_{-v} = \infty] = 0.$$

It thus follows from (2.39) and (2.41) that $P_0(A_l) > 0$. □

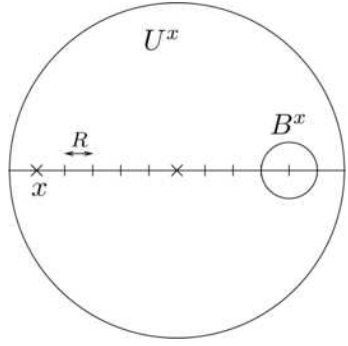

Fig. 2.

COROLLARY 2.6 ($d \geq 1$). *Let $l \in S^{d-1}$. If $P_0[A_l \cup A_{-l}] = 0$, then $P_0$-a.s., $\lim_{t \to \infty} \frac{l \cdot X_t}{t} = 0$.*

PROOF. If $P_0[A_l] = 0$, Theorem 2.4 implies that $\limsup_t \frac{l \cdot X_t}{t} \leq 0$, $P_0$-a.s. The same argument for $-l$ implies that $\liminf_t \frac{l \cdot X_t}{t} \geq 0, P_0$-a.s., and the claim follows. □

**3. Limit velocity.** The aim of this section is to prove the existence of a possibly nondeterministic asymptotic velocity; see Theorem 3.8. As a preparation we need to revisit some of the theorems proven in [26], now in the absence of the assumption $P_0(A_l) = 1$ made in [26]. For this reason we will consider in the following probabilities conditioned on the event that the diffusion is unbounded in a direction $l$ or that it escapes to infinity in a direction $l$. But first we recall the definitions of the regeneration times $\tau_k, k \geq 1$, and the coupling measure $\hat{P}_{x,\omega}$ introduced in [26].

3.1. *The coupling measure.* For $x \in \mathbb{R}^d$ and $l \in S^{d-1}$, we consider

(3.1) $\qquad B^x = B(x + 9Rl, R), \qquad U^x = B(x + 5Rl, 6R),$

where $R$ is the range of dependence of the environment. See Figure 2.

We denote by $\lambda_j$ the canonical coordinates on $\{0,1\}^{\mathbb{N}}$. Further, we let $(\mathcal{S}_m)_{m \geq 0}$ denote the canonical filtration on $\{0,1\}^{\mathbb{N}}$ and $\mathcal{S}$ the canonical $\sigma$-algebra. On the enlarged space $C(\mathbb{R}_+, \mathbb{R}^d) \times \{0,1\}^{\mathbb{N}}$, we consider the following $\sigma$-fields:

(3.2) $\quad \mathcal{Z}_t \stackrel{\text{def}}{=} \mathcal{F}_t \otimes \mathcal{S}_{\lceil t \rceil} \qquad \text{with } t \geq 0, \quad \text{and} \quad \mathcal{Z} \stackrel{\text{def}}{=} \mathcal{F} \otimes \mathcal{S} = \sigma\left\{ \bigcup_{m \in \mathbb{N}} \mathcal{Z}_m \right\}.$

On the enlarged space, the shift operators $\hat{\theta}_m$, $m \geq 0$, are defined so that $\hat{\theta}_m(X., \lambda.) = (X_{m+\cdot}, \lambda_{m+\cdot})$. Then from Theorem 2.1 in [26], one has the following measures, coupling the diffusion in random environment with a sequence of Bernoulli variables:

PROPOSITION 3.1. *There exists $\varepsilon > 0$, such that for every $l \in S^{d-1}$, $\omega \in \Omega$ and $x \in \mathbb{R}^d$, there exists a probability measure $\hat{P}_{x,\omega}$ on $(C(\mathbb{R}_+, \mathbb{R}^d) \times \{0,1\}^{\mathbb{N}}, \mathcal{Z})$ depending measurably on $\omega$ and $x$, such that:*

(3.3) *Under $\hat{P}_{x,\omega}$, $(X_t)_{t \geq 0}$ is $P_{x,\omega}$-distributed, and the $\lambda_m$, $m \geq 0$, are i.i.d. Bernoulli variables with success probability $\varepsilon$.*

(3.4) *For $m \geq 1$, $\lambda_m$ is independent of $\mathcal{F}_m \otimes \mathcal{S}_{m-1}$ under $\hat{P}_{x,\omega}$. Conditioned on $\mathcal{Z}_m$, $X_\cdot \circ \hat{\theta}_m$ has the same law as $X_\cdot$ under $\hat{P}^{\lambda_m}_{X_m,\omega}$, where for $y \in \mathbb{R}^d$, $\lambda \in \{0,1\}$, $\hat{P}^\lambda_{y,\omega}$ denotes the law $\hat{P}_{y,\omega}[\cdot | \lambda_0 = \lambda]$.*

(3.5) *$\hat{P}^1_{x,\omega}$ almost surely, $X_s \in U^x$ for $s \in [0,1]$ [recall (3.1)].*

(3.6) *Under $\hat{P}^1_{x,\omega}$, $X_1$ is uniformly distributed on $B^x$ [recall (3.1)].*

We then introduce the new annealed measures on $(\Omega \times C(\mathbb{R}_+, \mathbb{R}^d) \times \{0,1\}^{\mathbb{N}}, \mathcal{A} \otimes \mathcal{Z})$:

$$\hat{P}_x \stackrel{\text{def}}{=} \mathbb{P} \times \hat{P}_{x,\omega} \quad \text{and} \quad \hat{E}_x \stackrel{\text{def}}{=} \mathbb{E} \times \hat{E}_{x,\omega}. \tag{3.7}$$

3.2. *The regeneration times $\tau_k$.* We follow [26] and [25] to define the first regeneration time $\tau_1$. To this end, we introduce a sequence of integer-valued $(\mathcal{Z}_t)_{t \geq 0}$-stopping times $N_k$, for which the condition $\lambda_{N_k} = 1$ holds, and at these times the process $(l \cdot X_s)_{s \geq 0}$ reaches essentially a local maximum (within a small variation). Then $\tau_1$, when finite, is the first $N_k + 1$, $k \geq 1$, such that the process $(l \cdot X_t)_{t \geq 0}$ never goes below $l \cdot X_{N_k+1} - R$ after time $N_k + 1$. In fact, the precise definition of $\tau_1$ relies on several sequences of stopping times. First, for $a > 0$, introduce the $(\mathcal{F}_t)_{t \geq 0}$-stopping times $V_k(a)$, $k \geq 0$ [recall $T_u$ in (2.2)]:

$$V_0(a) \stackrel{\text{def}}{=} T_{M(0)+a} \leq \infty, \qquad V_{k+1}(a) \stackrel{\text{def}}{=} T_{M(\lceil V_k(a) \rceil) + R} \leq \infty \tag{3.8}$$

where $M(t) \stackrel{\text{def}}{=} \sup\{l \cdot X_s : 0 \leq s \leq t\}$.

In view of the Markov property [see (3.4)], we require the stopping times $N_k(a)$, $k \geq 1$, to be integer-valued and with this in mind, introduce as an intermediate step the (integer-valued) stopping times $\tilde{N}_k(a)$ where the process $X_t \cdot l$ essentially reaches a maximum:

$$\tilde{N}_1(a) \stackrel{\text{def}}{=} \inf\left\{\lceil V_k(a) \rceil : k \geq 0, \sup_{s \in [V_k, \lceil V_k \rceil]} |l \cdot (X_s - X_{V_k})| < \frac{R}{2}\right\},$$

(3.9)
$$\tilde{N}_{k+1}(a) \stackrel{\text{def}}{=} \tilde{N}_1(3R) \circ \hat{\theta}_{\tilde{N}_k(a)} + \tilde{N}_k(a), \qquad k \geq 1.$$

By convention we set $\tilde{N}_0 = 0$ and $\tilde{N}_{k+1} = \infty$ if $\tilde{N}_k = \infty$ and then define $N_1(a)$ as

$$N_1(a) \stackrel{\text{def}}{=} \inf\{\tilde{N}_k(a) : k \geq 1, \lambda_{\tilde{N}_k(a)} = 1\}. \tag{3.10}$$

Now we can define the $(\mathcal{Z}_t)_{t \geq 0}$-stopping times:

$$S_1 \stackrel{\text{def}}{=} N_1(3R) + 1, \qquad R_1 \stackrel{\text{def}}{=} S_1 + D \circ \hat{\theta}_{S_1} \tag{3.11}$$

$$\text{with } D \stackrel{\text{def}}{=} \lceil \tilde{T}^{\text{rel}}_{-R} \rceil. \tag{3.12}$$

(By convention we set $R_0 = 0$.) The $(\mathcal{Z}_t)_{t \geq 0}$-stopping times $N_{k+1}$, $S_{k+1}$ and $R_{k+1}$ are defined in an iterative way for $k \geq 1$:

$$\begin{aligned} N_{k+1} &\stackrel{\text{def}}{=} R_k + N_1(a_k) \circ \hat{\theta}_{R_k} \qquad \text{with } a_k \stackrel{\text{def}}{=} M(R_k) - X_{R_k} \cdot l + R, \\ S_{k+1} &\stackrel{\text{def}}{=} N_{k+1} + 1, \qquad R_{k+1} \stackrel{\text{def}}{=} S_{k+1} + D \circ \hat{\theta}_{S_{k+1}} \end{aligned} \tag{3.13}$$

(the shift $\hat{\theta}_{R_k}$ is *not* applied to $a_k$ in the above definition).

For $k \geq 1$, observe that on the event $\{N_k < \infty\}$, $\lambda_{N_k} = 1$ and $\sup_{s \leq N_k} X_s \cdot l \leq X_{N_k} \cdot l + R$. Notice that for all $k \geq 1$, the $(\mathcal{Z}_t)_{t \geq 0}$-stopping times $N_k$, $S_k$ and $R_k$ are integer-valued, possibly equal to infinity, and we have $1 \leq N_1 \leq S_1 \leq R_1 \leq N_2 \leq S_2 \leq R_2 \leq \cdots \leq \infty$.

The first *regeneration time* $\tau_1$ is defined, as in [25] and [26] (see also [33]), by

$$\tau_1 \stackrel{\text{def}}{=} \inf\{S_k : S_k < \infty, R_k = \infty\} \leq \infty. \tag{3.14}$$

3.3. *Renewal structure and limit velocity.* We first develop the main theorems describing the renewal structure and then present a weak zero–one law, which says that for any unit vector $l$, $P_0(A_l \cup A_{-l})$ is either 0 or 1; see Proposition 3.6. We then prove finiteness of $\hat{E}_0[l \cdot X_{\tau_1} | D = \infty]$ under the condition $P_0(A_l) > 0$ (cf. Proposition 3.7) and derive the existence of a possibly random asymptotic velocity in Theorem 3.8. We begin with an easy lemma which refines (2.41).

LEMMA 3.2. *For any $l \in S^{d-1}$, $P_0(A_l) > 0 \Leftrightarrow P_0(\tilde{T}_{-R} = \infty) > 0$.*

PROOF. In view of (2.41), we only need to prove that $P_0(A_l) > 0$ implies $P_0(\tilde{T}_{-R} = \infty) > 0$. Assume by contradiction that $P_0(\tilde{T}^{\text{rel}}_{-R} = \infty) = 0$. Using translation invariance of $\mathbb{P}$, and Fubini's theorem, we see that for almost all $\omega \in \Omega$

$$P_{x,\omega}(\tilde{T}^{\text{rel}}_{-R} = \infty) = 0 \qquad \text{except on a Lebesgue-negligible subset of } \mathbb{R}^d.$$

A calculation similar to (2.38) shows that for almost all $\omega$ and every $x \in \mathbb{R}^d$, we have that $P_{x,\omega}(\tilde{T}^{\text{rel}}_{-R/2} < \infty) = 1$. The strong Markov property implies at once that $P_{0,\omega}(\tilde{T}^{\text{rel}}_{-kR/2} < \infty) = 1$, $\mathbb{P}$-a.s., for all $k \geq 1$. This contradicts $P_0(A_l) > 0$. $\square$

In the sequel, we will use the following additional notation. For any $l \in S^{d-1}$,

$$B_l \stackrel{\text{def}}{=} \left\{ \sup_{s \geq 0} l \cdot X_s = \infty \right\}. \tag{3.15}$$

From (2.41) and the definition of $D$ [see (3.12)], we have of course for any $l \in S^{d-1}$

$$\{D = \infty\} \subset A_l \subset B_l, \qquad P_0\text{-a.s.} \tag{3.16}$$

We will see later that if $P_0(A_l) > 0$, then $A_l = B_l$, $P_0$-a.s.; see Theorem 3.5. The next lemma shows that the first renewal time $\tau_1$ is finite on the event $B_l$, if $P_0(A_l) > 0$.

LEMMA 3.3. *Consider $l \in S^{d-1}$ and assume $P_0(A_l) > 0$; then $B_l \subset \{\tau_1 < \infty\}$, $\hat{P}_0$-a.s., with the notation* (3.15).

The proof is similar to the proof of Proposition 2.7 in [26] and is included in the Appendix for the convenience of the reader.

On the space $\Omega \times C(\mathbb{R}_+, \mathbb{R}^d) \times \{0,1\}^{\mathbb{N}}$, we introduce the sub-$\sigma$-algebra $\mathcal{G}$ of $\mathcal{A} \otimes \mathcal{Z}_\infty$ that is generated by sets of the form

$$\{\tau_1 = m\} \cap O_{m-1} \cap \{l \cdot X_{m-1} > a\} \cap \{X_m \in G\} \cap F_a, \tag{3.17}$$

where $m \geq 2, a \in \mathbb{R}, O_{m-1} \in \mathcal{Z}_{m-1}, G \subset \mathbb{R}^d$ open, $F_a \in \mathcal{H}_{\{z \in \mathbb{R}^d : l \cdot z \leq a+R\}}$.

Loosely speaking, $\mathcal{G}$ contains information on the trajectories up to time $\tau_1 - 1$ and at time $\tau_1$ as well as information on the environment that has possibly been visited by the diffusion up to time $\tau_1 - 1$. Note that no information between time $\tau_1 - 1$ and $\tau_1$ is included. This is crucial when one exploits the finite range dependence property of the environment with the help of the coupling measure $\hat{P}_0$, as we saw already in the proof of Lemma 3.3; see (A.2).

The next proposition is a variation on Theorem 2.4 in [26], and provides the base for the renewal structure presented in Theorem 3.5.

PROPOSITION 3.4. *Consider $l$ in $S^{d-1}$ and assume $P_0(A_l) > 0$. Then for any $x \in \mathbb{R}^d$, any bounded functions $f, g, h$ respectively $\mathcal{Z}, \mathcal{H}_{\{z \in \mathbb{R}^d : l \cdot z \geq -R\}}$, $\mathcal{G}$-measurable, one has*

$$\begin{aligned}
\hat{E}_x[f(X_{\tau_1 + \cdot} - X_{\tau_1}, \lambda_{\tau_1 + \cdot}) g \circ t_{X_{\tau_1}} h | B_l] \\
= \hat{E}_x[h | B_l] \hat{E}_0[f(X_\cdot, \lambda_\cdot) g | D = \infty],
\end{aligned} \tag{3.18}$$

with $B_l$ as in (3.15) and $t_y$ the spatial shift; see the beginning of the Introduction.

[We will later see that $A_l = B_l$, if $P_0(A_l) > 0$; see Theorem 3.5.]

PROOF. We only discuss the salient features of the proof, which is a variation on that of Theorem 2.4 in [26]. As in the proof of this theorem, it suffices to prove (3.18) for $h = \mathbb{1}_{\{\tau_1 = m\}} \mathbb{1}_{F_a} \mathbb{1}_{O_{m-1}} \mathbb{1}_{\{X_m \in G\}} \mathbb{1}_{\{l \cdot X_{m-1} > a\}}$, with $m \geq 2, a \in \mathbb{R}, O_{m-1} \in \mathcal{Z}_{m-1}, G \subset \mathbb{R}^d$ open, $F_a \in \mathcal{H}_{\{z \in \mathbb{R}^d : l \cdot z \leq a + R\}}$, since (3.17) constitutes a $\pi$-system. Note that there is a $\tilde{O}_{m-1} \in \mathcal{Z}_{m-1}$, such that

$$\{\tau_1 = m\} \cap O_{m-1} \cap B_l = \tilde{O}_{m-1} \cap \{D \circ \hat{\theta}_m = \infty\} \cap \hat{\theta}_m^{-1}(B_l)$$
$$= \tilde{O}_{m-1} \cap \hat{\theta}_m^{-1}(\{D = \infty\} \cap B_l)$$
$$= \tilde{O}_{m-1} \cap \hat{\theta}_m^{-1}(\{D = \infty\}),$$

where the last step follows from (3.16). As $B_l$ disappears from the calculations, the rest of the argument is identical to the proof of Theorem 2.4 in [26]; see also [27]. □

On the event $\{\tau_1 < \infty\}$, we define inductively a nondecreasing sequence of random variables $\tau_k \leq \infty$, via

$$(3.19) \quad \tau_{k+1}((X., \lambda.)) \stackrel{\text{def}}{=} \tau_1((X., \lambda.)) + \tau_k((X_{\tau_1 + \cdot}, \lambda_{\tau_1 + \cdot})), \qquad k \geq 1.$$

We are able to reconstruct in our context an analogue of the renewal structure of Theorem 2.5 in [26].

THEOREM 3.5. *Consider $l$ in $S^{d-1}$ and assume $P_0(A_l) > 0$. Then $\hat{P}_0$-a.s., $\{D = \infty\} \subset A_l = B_l = \{\tau_k < \infty, \text{ for all } k \geq 1\}$ [recall (3.12), (1.9), (3.16)] and under $\hat{P}_0[\cdot | A_l]$, the random variables*

$$(3.20) \quad Z_k \stackrel{\text{def}}{=} (X_{(\tau_k + \cdot) \wedge (\tau_{k+1} - 1)} - X_{\tau_k}, X_{\tau_{k+1}} - X_{\tau_k}, \tau_{k+1} - \tau_k), \qquad k \geq 0,$$

*are independent. Moreover under $\hat{P}_0[\cdot | A_l]$, the random variables $Z_k$, $k \geq 1$ have the same distribution as $Z_0$ under $\hat{P}_0[\cdot | D = \infty]$.*

PROOF. We use induction over the index $n \geq 0$, of the filtration $\mathcal{G}_n \stackrel{\text{def}}{=} \sigma(Z_0, \ldots, Z_n)$. From Proposition 3.4 and the fact that $\mathcal{G}_0 \subset \mathcal{G}$ [cf. (3.17)], we know that for any $C$ in the product $\sigma$-algebra on $C(\mathbb{R}_+, \mathbb{R}^d) \times \mathbb{R}^d \times \mathbb{R}_+$ and any bounded $\mathcal{G}_0$-measurable $h_0$,

$$(3.21) \qquad \hat{E}_0[\mathbb{1}_{\{Z_1 \in C\}} h_0 | B_l] = \hat{E}_0[h_0 | B_l] \hat{P}_0[Z_0 \in C | D = \infty].$$

It follows that on $B_l$, $\tau_2$ is $\hat{P}_0$-a.s. finite, because $\hat{P}_0[\tau_2 < \infty | B_l] = \hat{P}_0[\tau_1 < \infty | D = \infty] = 1$ by Lemma 3.3 and (3.16). Assume now that for some $n \geq$

1, $\tau_n < \infty$, on $B_l$ and that for any $C$ as above and any bounded $\mathcal{G}_{n-1}$-measurable $h_{n-1}$:

$$(3.22) \quad \hat{E}_0[\mathbb{1}_{\{Z_n \in C\}} h_{n-1} | B_l] = \hat{E}_0[h_{n-1} | B_l] \hat{P}_0[Z_0 \in C | D = \infty].$$

As above we see that $\hat{P}_0[\tau_{n+1} < \infty | B_l] = 1$. We will prove an identity similar to (3.22) with $(n+1)$ in place of $n$. By the definition of $\tau_{n+1}$, $\mathcal{G}_n \cap \{\tau_1 < \infty\}$ is generated by a $\pi$-system consisting of intersections between events in $\mathcal{G}_0 \cap \{\tau_1 < \infty\}$ and $\hat{\theta}_{\tau_1}^{-1} \mathcal{G}_{n-1}$. With Dynkin's lemma (see [7], page 447), it suffices therefore to consider bounded, $\mathcal{G}_n$-measurable functions $h_n$ satisfying

$$(3.23) \quad h_n = h_0 \cdot h_{n-1} \circ \hat{\theta}_{\tau_1}, \qquad \hat{P}_0\text{-a.s. on } \{\tau_1 < \infty\},$$

for some bounded $\mathcal{G}_0$-measurable, respectively $\mathcal{G}_{n-1}$-measurable, functions $h_0$ and $h_{n-1}$. Let us now prove the induction step with $h_n$ as in (3.23). By Proposition 3.4, we have for any $C$ as above

$$(3.24) \quad \begin{aligned} \hat{E}_0[\mathbb{1}_{\{Z_{n+1} \in C\}} h_n | B_l] &= \hat{E}_0[(h_{n-1} \mathbb{1}_{\{Z_n \in C\}}) \circ \hat{\theta}_{\tau_1} h_0 | B_l] \\ &= \frac{\hat{E}_0[h_0 | B_l]}{P_0[D = \infty]} \hat{E}_0[\mathbb{1}_{\{Z_n \in C\}} h_{n-1} \mathbb{1}_{\{D = \infty\}}]. \end{aligned}$$

Let us admit for the time being that

(3.25) $h_{n-1} \mathbb{1}_{\{D = \infty\}}$ is indistinguishable from a $\mathcal{G}_{n-1}$-measurable variable

and conclude the induction step. It follows from (3.22), (3.25) and the fact $P_0$-a.s., $\{D = \infty\} \subset B_l$ [cf. (3.16)] that the left-hand side of (3.24) equals

$$(3.26) \quad \hat{E}_0[h_0 | B_l] \hat{E}_0[h_{n-1} | D = \infty] \hat{P}_0[Z_0 \in C | D = \infty].$$

Replacing $C$ with $C(\mathbb{R}_+, \mathbb{R}^d) \times \mathbb{R}^d \times \mathbb{R}$, we obtain

$$\hat{E}_0[h_n | B_l] = \hat{E}_0[h_0 | B_l] \hat{E}_0[h_{n-1} | D = \infty].$$

Inserting this into (3.26) yields $\hat{E}_0[\mathbb{1}_{\{Z_{n+1} \in C\}} h_n | B_l] = \hat{E}_0[h_n | B_l] \hat{P}_0[Z_0 \in C | D = \infty]$. In other words, (3.22) holds with $(n+1)$ in place of $n$. Note that the induction argument shows that if $P_0(A_l) > 0$, then $\hat{P}_0$-a.s., $B_l \subset \{\tau_k < \infty$, for all $k \geq 0\}$ and thus $\hat{P}_0$-a.s., $B_l = A_l = \{\tau_k < \infty$, for all $k \geq 0\}$. [We will see later that in fact $\hat{P}_0$-a.s, $A_l = \{\tau_1 < \infty\}$, if $P_0(A_l) > 0$; cf. Proposition 3.6.]

It remains to prove (3.25): Observe that $\hat{P}_0$-a.s., $\{D = \infty\} = \{\tilde{T}_{-R} = \infty\} = \{\tilde{T}_{-R} \geq \tau_1\}$. Further it is clear that the last event is included in $\{D \geq \tau_1\}$; see (3.12). They are in fact equal $\hat{P}_0$-a.s., because the converse inclusions stems from the following facts: $\{D \geq \tau_1\} \cap \{\tau_1 = \infty\}$ is a $\hat{P}_0$ null-set by (3.16) and Lemma 3.3, and $\{D \geq \tau_1\} \cap \{\tau_1 < \infty\} \subset \{\tilde{T}_{-R} > \tau_1 - 1\}, \hat{P}_0$-a.s. But on $\{\tau_1 < \infty\}, \hat{P}_0$-a.s. $l \cdot X_{\tau_1 - 1 + s} \geq 2R$, for all $s \geq 0$, by construction of $\tau_1$. Therefore $\{\tilde{T}_{-R} > \tau_1 - 1\} \subset \{\tilde{T}_{-R} = \infty\} \subset \{\tilde{T}_{-R} \geq \tau_1\}$, $\hat{P}_0$-a.s.

We thus see that $\hat{P}_0$-a.s., $\{D = \infty\} = \{D \geq \tau_1\} = \{D \leq \tau_1 - 1\}^c$ which is $\mathcal{G}_0$-measurable and thus $h_{n-1}\mathbb{1}_{\{D=\infty\}}$ is indistinguishable from a $\mathcal{G}_{n-1}$-measurable variable. $\square$

PROPOSITION 3.6 (Weak zero–one law, $d \geq 1$). *For any $l \in S^{d-1}, P_0(A_l \cup A_{-l}) \in \{0,1\}$. Moreover if $P_0(A_l) > 0$, then $\hat{P}_0$-a.s., $B_l = A_l = \{\tau_1 < \infty\}$, where $B_l$ is defined in* (3.15).

PROOF. Assume that $P_0(A_l) > 0$, and consider any $L > 0$. Let $H_k, k \geq 0$, be the iterates of $H_{\bar{\mathcal{S}}_{(-L,L)}} \circ \theta_1 + 1$. We claim that $P_0[H_k < \infty, \text{for all } k \geq 0] = 0$. Indeed, using the notation from Lemma 2.5, we see that

$$P_0[\{H_k < \infty, \text{for all } k \geq 0\} \cap B_l^c]$$
$$\leq P_0\left[\bigcup_{v \in \mathbb{N}} \{\beta_k^{l,-L} < \infty, \text{ for all } k \geq 0 \text{ and } T_v = \infty\}\right] = 0.$$

From Theorem 3.5, we know that $B_l = A_l, P_0$-a.s. and therefore we find

$$P_0[\{H_k < \infty, \text{for all } k \geq 0\} \cap B_l] = 0.$$

This proves the claim and as $L$ is arbitrary, we see that $P_0[\lim_{t \to \infty} |l \cdot X_t| = \infty] = 1$, and hence $P_0(A_l \cup A_{-l}) = 1$, under the assumption $P_0(A_l) > 0$. The case where $P_0(A_{-l}) > 0$ is treated analogously and the 0–1 law follows. Finally observe that under the assumption $P_0(A_l) > 0$, we have that $\hat{P}_0$-a.s., $\{\tau_1 < \infty\} = (\{\tau_1 < \infty\} \cap A_l) \cup (\{\tau_1 < \infty\} \cap A_{-l})$, where the second set in the union is empty. Hence $\{\tau_1 < \infty\} \subset A_l, \hat{P}_0$-a.s. The converse inclusion follows from Lemma 3.3. $\square$

The next proposition proves that $P_0(A_l) > 0$ implies that $l \cdot X_{\tau_1}$ has a finite first moment under $\hat{P}_0[\cdot|D = \infty]$. In the discrete i.i.d. setting where the renewal structure is technically less intricate (cf., e.g., [33]) and under the assumption that $l$ is a coordinate direction, one can show a stronger result, namely the equality $E_0[l \cdot X_{\tau_1}|D = \infty] = P_0[D = \infty]^{-1}$; see, for instance, [36], Lemma 3.2.5. Let us now give an outline of the argument we use. We find an $L > 0$ for which there is a positive lower bound, uniform in $r > 0$, for the annealed probability that the interval $[r, r+L]$ contains one of the $l \cdot X_{\tau_m}, m \geq 1$. This yields a positive lower bound on the linear growth in $M$ of the expected number of renewal points $l \cdot X_{\tau_m}$ smaller than $M$. But by the elementary renewal theorem this linear growth coincides with $\hat{E}_0[l \cdot X_{\tau_1}|D = \infty]^{-1}$. We thus obtain the desired upper bound on $\hat{E}_0[l \cdot X_{\tau_1}|D = \infty]$. Whereas the construction of an $L$ as above is relatively straightforward in the discrete setup, it is somewhat involved in the continuous setting because of the more delicate nature of the regeneration times.

Let us incidentally point out that the use of the elementary renewal theorem bypasses the arithmeticity conditions of Blackwell's renewal theorem used in [36]. This is an advantage when working with a general direction $l$ (both in the discrete and continuous setups).

PROPOSITION 3.7. *Consider $l$ in $S^{d-1}$ and assume that $P_0(A_l) > 0$. Then there is a constant $c_0 > 0$ such that if $L$ is large enough, for any $r \geq 0$ one has*

(3.27) $\qquad \hat{P}_0[\text{ for some } m \geq 1, l \cdot X_{\tau_m} \in [r, r+L] | A_l] > c_0,$

(3.28) $$\hat{E}_0[l \cdot X_{\tau_1} | D = \infty] \leq \frac{L}{c_0}.$$

PROOF. We first prove (3.27). Consider any $r \geq 0$, $0 < \delta < \frac{R}{10}$ and define $T = T_{r+R/4}$; see (2.2). The heart of the matter is to construct an event $E$ [cf. (3.30)], forcing the occurrence of some $l \cdot X_{\tau_m}, m \geq 1$, in an interval. More precisely we will show that

(3.29) $\hat{P}_0$-a.s., on the event $E$, some $l \cdot X_{\tau_m}, m \geq 1$, belongs to $[r, r+18R]$,

where $E$ is defined as

$$E = \Big\{ T < \infty, \sup_{s \in [T, \lceil T \rceil]} |X_s - X_T| \leq \frac{R}{4},$$

(3.30) $$\sup_{0 \leq s \leq 2} |X_s - X_0 - \psi(s)| \circ \hat{\theta}_{\lceil T \rceil} < \delta,$$

$$\lambda_{\lceil T \rceil + 2} = 1, D \circ \hat{\theta}_{\lceil T \rceil + 3} = \infty \Big\},$$

and $\psi : \mathbb{R}_+ \to \mathbb{R}^d$ is the function

$$\psi(s) = \begin{cases} 5Rls, & s \leq 1, \\ 5Rl + (s-1)\frac{5}{4}Rl, & 1 < s \leq 2. \end{cases}$$

The intuitive idea behind the construction of $E$ is the following (see Figure 3): in essence after first reaching level $r + R/4$ at time $T$, the trajectory is forced—in the next unit of time after $\lceil T \rceil$—to move $5R$ "to the right" and—in the subsequent unit of time—to move an additional distance $R$ "to the right." Then either $\lceil T \rceil$ coincides with a regeneration time, or as we will see, some time "of type V" [after suitable time shift, see (3.8)] occurs during the first interval $[\lceil T \rceil, \lceil T \rceil + 1]$. Rounding up this time to the next integer yields $\lceil T \rceil + 1$ and the constraints imposed on the trajectory during the second unit of time $[\lceil T \rceil + 1, \lceil T \rceil + 2]$ as well as on the Bernoulli variables, ensure that $\lceil T \rceil + 2$ is "of type N"; see (3.10). Because of the no-backtracking condition in $E$, $\lceil T \rceil + 3$ is then a regeneration time and we have a good control on how far "to the right" the trajectory has moved at that time.

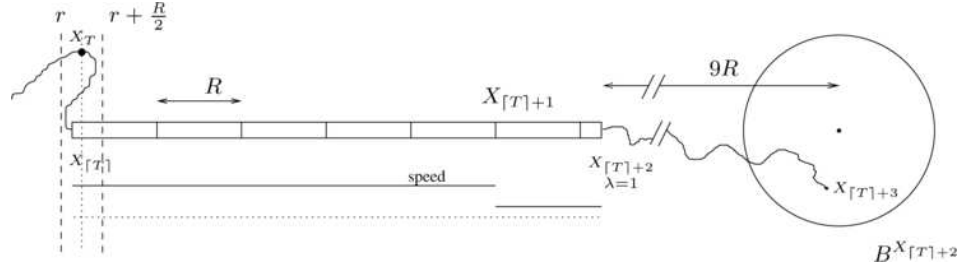

FIG. 3. *A realization of the event $E$ [cf. (3.30)] and the corresponding speed of the trajectory.*

We now proceed with the proof of (3.29). Let $\tau_m < \lceil T \rceil, m \geq 0$, be the last regeneration time strictly before $\lceil T \rceil$, with $m = 0$ by convention when $\lceil T \rceil = 0$, which is a $\hat{P}_0$-negligible event. We define

$$(3.31) \qquad k = \sup\{n \geq 0 : R_n \circ \hat{\theta}_{\tau_m} + \tau_m \leq \lceil T \rceil\};$$

see (3.11), (3.13) for the notation. On the event $E$, the following two cases can occur:

(i) Either $N_{k+1} \circ \hat{\theta}_{\tau_m} + \tau_m < \lceil T \rceil$, then we claim that

$$(3.32) \qquad N_{k+1} \circ \hat{\theta}_{\tau_m} + \tau_m + 1 = \tau_{m+1} = \lceil T \rceil.$$

Indeed, according to the definition (3.14), (3.19) of $\tau_{m+1}$, the first equality in (3.32) automatically holds if $R_{k+1} \circ \hat{\theta}_{\tau_m}$ is infinite. Assume by contradiction, that $R_{k+1} \circ \hat{\theta}_{\tau_m} < \infty$. By the definition of $k$, $R_{k+1} \circ \hat{\theta}_{\tau_m} + \tau_m > \lceil T \rceil$ and by the definition of $R_{k+1}$ [cf. (3.11)], the trajectory would have to return to level $u^* \stackrel{\text{def}}{=} l \cdot X_{N_{k+1} \circ \hat{\theta}_{\tau_m} + \tau_m + 1} - R$ strictly after time $\lceil T \rceil$. But under our assumption, $N_{k+1} \circ \hat{\theta}_{\tau_m} + \tau_m + 1 \leq \lceil T \rceil$ and hence with the second condition in the definition of $E$, $u^* \leq l \cdot X_{\lceil T \rceil} - \frac{3R}{4}$. On $E$, however, after time $\lceil T \rceil$, the trajectory always stays strictly above level $l \cdot X_{\lceil T \rceil} - \frac{R}{2}$. This contradiction proves that $R_{k+1} \circ \hat{\theta}_{\tau_m}$ is infinite and hence the first equality of (3.32) follows. The second equality simply stems from the fact that $N_{k+1} \circ \hat{\theta}_{\tau_m} + \tau_m + 1 \leq \lceil T \rceil \leq \tau_{m+1}$ in the considered case.

(ii) Or $\lceil T \rceil \leq N_{k+1} \circ \hat{\theta}_{\tau_m} + \tau_m$; then we first note that $\lceil T \rceil = N_{k+1} \circ \hat{\theta}_{\tau_m} + \tau_m$ is $\hat{P}_0$-negligible as $\{\lambda_{\lceil T \rceil} = 1\} \cap E$ is a $\hat{P}_0$ null-set by (3.6). We claim that

$$(3.33) \qquad \lceil T \rceil + 3 = \tau_{m+1}.$$

To see this, we first determine below a random time $\bar{N} \leq \lceil T \rceil$ "of type $\tau$, $R$ or $\tilde{N}$," serving as starting point for the construction of a new generation of stopping times "of type $V$;" see Section 3.2. With $k$ as in (3.31), we define

$$\rho = R_k \circ \hat{\theta}_{\tau_m} + \tau_m, \qquad \bar{N}_j = \tilde{N}_j \circ \hat{\theta}_\rho + \rho, j \geq 0,$$
$$j_0 = \sup\{j \geq 0 : \bar{N}_j \leq \lceil T \rceil\}, \qquad \bar{N} = \bar{N}_{j_0}.$$

Since on $E$ the trajectory visits a new half plane once it reaches level $r + \frac{R}{2}$, there exists a smallest $i \geq 0$, such that $V \stackrel{\text{def}}{=} V_i(a) \circ \hat{\theta}_{\bar{N}} + \bar{N}$ [where $a$ equals either $3R$ or $M(\rho) - X_\rho \cdot l + R$ according to the type of $\bar{N}$; cf. (3.8), (3.9), (3.11), (3.13)], satisfies

$$\lceil T \rceil < V \leq (N_{k+1} \circ \hat{\theta}_{\tau_m} + \tau_m) \wedge (\lceil T \rceil + 1)$$

and

$$l \cdot X_{\lceil T \rceil} < l \cdot X_V < l \cdot X_{\lceil T \rceil} + 4R.$$

Note that $\lceil V \rceil = \lceil T \rceil + 1 \neq \bar{N}_{j_0+1}$ because $l \cdot (X_{\lceil V \rceil} - X_V) > \frac{R}{2}$; see also (3.9). But the "next" $V$, namely $V' \stackrel{\text{def}}{=} V_{i+1}(a) \circ \hat{\theta}_{\bar{N}} + \bar{N} > \lceil T \rceil + 1$, is reached by definition of $E$ near level $l \cdot X_{\lceil T \rceil} + 6R$, and $\lceil V' \rceil$ coincides with

$$N_{k+1} \circ \hat{\theta}_{\tau_m} + \tau_m = \lceil T \rceil + 2,$$

since $l \cdot (X_{\lceil V' \rceil} - X_{V'}) \leq \frac{R}{2}$ and $\lambda_{\lceil T \rceil + 2} = 1$. We obtain that $\lceil T \rceil + 3$ is the next regeneration time $\tau_{m+1}$, since on $E$ the trajectory never backtracks after $\lceil T \rceil + 3$. This proves (3.33).

So far we have shown (3.29) and there remains to prove that the probability $\hat{P}_0[E|A_l]$ is bounded away from 0, independently of $r$. The claim (3.27) will then follow. To this end, we observe that

$$\hat{P}_0[E \cap A_l] = \sum_{n=0}^{\infty} \mathbb{E}\hat{E}_{0,\omega}\bigg[\lceil T \rceil = n, \sup_{T \leq s \leq \lceil T \rceil} |l \cdot (X_s - X_T)| \leq \frac{R}{4},$$

(3.34)
$$\sup_{0 \leq s \leq 2} |X_s - X_0 - \psi(s)| \circ \hat{\theta}_n < \delta, \lambda_{n+2} = 1,$$

$$\hat{P}_{0,\omega}[D \circ \hat{\theta}_{n+3} = \infty, \hat{\theta}_{n+3}^{-1}(A_l)|\mathcal{Z}_{n+2}]\bigg].$$

With the Markov property (3.4) as well as (3.6) and the first inclusion in (3.16), we find that for $\mathbb{P}$-a.e. $\omega$

$$\hat{P}_{0,\omega}[D \circ \hat{\theta}_{n+3} = \infty, \hat{\theta}_{n+3}^{-1}(A_l)|\mathcal{Z}_{n+2}]$$

(3.35)
$$= \hat{P}_{X_{n+2},\omega}^{\lambda_{n+2}}[D \circ \hat{\theta}_1 = \infty]$$

$$= \frac{1}{|B(0,R)|} \int P_{y,\omega}[D = \infty] \mathbb{1}_{\{B^{X_{n+2}}\}}(y) \, dy.$$

We insert (3.35) into (3.34) and use the following facts:

(a) $\{\lambda_{n+2} = 1\}$ has probability $\varepsilon$ and is independent of $\mathcal{F}_{n+2} \otimes \mathcal{S}_{n+1}$; see (3.3) and (3.4).

(b) $\omega \mapsto P_{0,\omega}[\lceil T \rceil = n, \sup_{T \leq s \leq \lceil T \rceil} |l \cdot (X_s - X_T)| \leq \frac{R}{4}, \sup_{0 \leq s \leq 2} |X_s - X_0 - \psi(s)| \circ \theta_n < \delta, y \in B^{X_{n+2}}]$ is $\mathcal{H}_{\{z \in \mathbb{R}^d : z \cdot l \leq y \cdot l - 8R + \delta\}}$-measurable [recall the definition of $B^x$, (3.1)]. $\omega \mapsto P_{y,\omega}[D = \infty]$ is $\mathcal{H}_{\{z \in \mathbb{R}^d : z \cdot l \geq y \cdot l - R\}}$-measurable. Therefore by the finite range dependence property (1.5) both maps are $\mathbb{P}$ independent.

Moreover from Lemma 3.2, we have that $P_0[D = \infty] > c > 0$ and hence we obtain

$$\hat{P}_0[E \cap A_l] = \frac{\varepsilon}{|B(0,R)|}$$

$$\times \sum_{n=0}^{\infty} \int_{\mathbb{R}^d} dy \, \mathbb{E} P_{0,\omega} \bigg[ \lceil T \rceil = n, \sup_{T \leq s \leq \lceil T \rceil} |l \cdot (X_s - X_T)| \leq \frac{R}{4},$$

$$\sup_{0 \leq s \leq 2} |X_s - X_0 - \psi(s)| \circ \theta_n < \delta, y \in B^{X_{n+2}} \bigg]$$

$$\times \mathbb{E} P_{y,\omega}[D = \infty]$$

$$\geq \varepsilon c \sum_{n=0}^{\infty} \mathbb{E} E_{0,\omega} \bigg[ \lceil T \rceil = n, \sup_{T \leq s \leq \lceil T \rceil} |l \cdot (X_s - X_T)| \leq \frac{R}{4},$$

$$P_{X_n,\omega} \bigg[ \sup_{0 \leq s \leq 2} |X_s - X_0 - \psi(s)| \leq \delta \bigg] \bigg]$$

$$\geq \varepsilon c c'(\delta, \psi) P_0[T < \infty] \geq \varepsilon c c'(\delta, \psi) P_0(A_l),$$

where the constant $c'(\delta, \psi)$ stems from the support theorem (see [2], page 25). This proves (3.27).

We now prove (3.28). From Theorem 3.5, we know that $l \cdot (X_{\tau_{m+1}} - X_{\tau_k})$, $m \geq 0$, are independent under $\hat{P}_0[\cdot | A_l]$ and have for $m \geq 1$ the same distribution as $l \cdot X_{\tau_1}$ under $\hat{P}_0[\cdot | D = \infty]$. Moreover $\hat{P}_0[\tau_1 < \infty | A_l] = 1$. Thus the elementary renewal theorem in the delayed case (see [24], Theorem 3.3.3) can be applied and yields

(3.36)
$$\hat{E}_0[l \cdot X_{\tau_1} | D = \infty]^{-1}$$
$$= \liminf_{k \to \infty} \frac{\hat{E}_0[\max\{m \geq 1 : l \cdot X_{\tau_m} \leq kL\} | A_l]}{kL} \geq \frac{c_0}{L}.$$

This proves (3.28). □

We now turn to the main result in this section, which describes the limiting velocity of the diffusion process.

THEOREM 3.8 (Limit velocity $d \geq 1$). *There exist a deterministic direction $l_* \in S^{d-1}$ and two numbers $v_+$, $v_- \geq 0$, such that*

(3.37) $\qquad P_0\text{-}a.s., \qquad \lim_{t \to \infty} \frac{X_t}{t} = (v_+ \mathbb{1}_{A_{l_*}} - v_- \mathbb{1}_{A_{-l_*}})l_*,$

*and $P(A_{l_*} \cup A_{-l_*}) \in \{0, 1\}$. (If this last quantity is 0, the velocity is 0 and thus the values of $v_+, v_-$ are immaterial.)*

PROOF. We first prove that for any fixed direction $l \in S^{d-1}$, there are nonnegative numbers $v_l, v_{-l}$, such that

(3.38) $\qquad P_0\text{-a.s.}, \qquad \lim_{t \to \infty} \frac{l \cdot X_t}{t} = v_l \mathbb{1}_{A_l} - v_{-l} \mathbb{1}_{A_{-l}}.$

If $P_0(A_l \cup A_{-l}) = 0$, it follows from Corollary 2.6 that (3.38) holds with $v_l = v_{-l} = 0$. In view of the weak zero–one law, Proposition 3.6, we only have to consider the case $P_0(A_l \cup A_{-l}) = 1$. We assume without loss of generality that $P(A_l) > 0$. On $A_l$, $\hat{P}_0$-a.s., $\tau_k < \infty$, $k \geq 1$ (cf. Theorem 3.5) and we define for $t > 0$, a nondecreasing, integer-valued function $k(t)$ tending to infinity $\hat{P}_0$-a.s., such that

$$\tau_{k(t)} \leq t < \tau_{k(t)+1},$$

with the convention $\tau_0 = 0$. Observe that on $A_l$, we have $\hat{P}_0$-a.s.,

(3.39) $\qquad \frac{l \cdot X_{\tau_{k(t)}} - R}{k(t)} \frac{k(t)}{\tau_{k(t)+1}} \leq \frac{l \cdot X_t}{t} \leq \frac{l \cdot X_{\tau_{k(t)+1}} + 3R}{k(t)+1} \frac{k(t)+1}{\tau_{k(t)}}.$

By (3.28), the i.i.d. structure of the increments $l \cdot (X_{\tau_k} - X_{\tau_{k-1}})$, $k \geq 2$, under $\hat{P}_0[\cdot | A_l]$ (see Theorem 3.5) and the usual law of large numbers, we find

(3.40) $\qquad \lim_{k \to \infty} \frac{l \cdot X_{\tau_k}}{k} = \hat{E}_0[l \cdot X_{\tau_1} | D = \infty] < \infty, \qquad \hat{P}_0[\cdot | A_l]\text{-a.s.}$

(i) Either $\hat{E}_0[\tau_1 | D = \infty] = \infty$, and then the positivity and the i.i.d. structure of the increments $\tau_k - \tau_{k-1}, k \geq 2$ (see Theorem 3.5) imply that $\frac{1}{n} \sum_{k=1}^{n} \tau_k - \tau_{k-1} \longrightarrow \infty, \hat{P}_0[\cdot | A_l]$-a.s. Passing to the limit in (3.39), we obtain in this case $v_l = 0$ in (3.38).

(ii) Or $\hat{E}_0[\tau_1 | D = \infty] < \infty$;

(3.41)
$$\text{then we obtain } v_l = \frac{\hat{E}_0[l \cdot X_{\tau_1} | D = \infty]}{\hat{E}_0[\tau_1 | D = \infty]} > 0.$$

If $P(A_{-l})$ is also positive, then the same argument determines $v_{-l}$; otherwise we set $v_{-l} = 0$. This proves (3.38).

Applying (3.38) to a basis of $\mathbb{R}^d$, we obtain

(3.42) $$X_t/t \longrightarrow v, \qquad P_0\text{-a.s.},$$

where $v$ is a random vector taking at most $2^d$ values.

In the next step we show that in fact $v$ takes at most two parallel and opposite values. Indeed, assume that there are $v_1, v_2$ noncolinear, nonzero vectors with $\hat{P}_0[v = v_i] > 0, i = 1, 2$. Define $e_i = \frac{v_i}{|v_i|}$, $i = 1, 2$ and

$$l_\alpha \stackrel{\text{def}}{=} \alpha e_1 + (1 - \alpha) e_2,$$

for $\alpha \in (0, 1)$. From (3.42) and (3.38) we see that $P_0$-a.s.,

$$v \cdot l_\alpha = v_{l_\alpha} \mathbb{1}_{A_{l_\alpha}} - v_{-l_\alpha} \mathbb{1}_{A_{-l_\alpha}} \qquad \text{for } \alpha \in (0, 1).$$

Therefore if for some $\alpha \in (0, 1)$, $l_\alpha \cdot v_i > 0$, for $i = 1, 2$, then since $\hat{P}_0[v = v_i] > 0$, we find

(3.43) $$l_\alpha \cdot v_1 = l_\alpha \cdot v_2.$$

If we can choose $\alpha$ in a nonempty open interval such that $l_\alpha \cdot v_i > 0, i = 1, 2$, holds, we may take derivatives with respect to $\alpha$ in (3.43) and deduce

$$0 = (e_1 - e_2)(v_1 - v_2)$$
$$= (e_1 - e_2)(|v_1|e_1 - |v_2|e_2)$$
$$= (1 - e_1 \cdot e_2)(|v_1| + |v_2|).$$

By assumption, $|e_1 \cdot e_2| < 1$, which produces a contradiction. Let us check that indeed $l_\alpha \cdot v_i > 0, i = 1, 2$, is true for $\alpha$ in a nonempty open interval:

$$l_\alpha \cdot v_1 > 0 \iff l_\alpha \cdot e_1 > 0 \iff \alpha > \frac{-e_1 \cdot e_2}{1 - e_1 \cdot e_2},$$
$$l_\alpha \cdot v_2 > 0 \iff l_\alpha \cdot e_2 > 0 \iff \alpha < \frac{1}{1 - e_1 \cdot e_2}.$$

Both bounds define a nonempty open interval as $|e_1 \cdot e_2| < 1$. As a result, there is an $l_* \in S^{d-1}$ such that $P_0[v \in \mathbb{R}l_*] = 1$. The application of (3.38) with $l_*$ together with (3.42) completes the proof. $\square$

**4. Zero–one law when $d = 2$.** In this section we prove that in the two-dimensional case, for any direction $l$, $P_0(A_l)$ is either 0 or 1. Note that this result combined with Theorem 3.8 implies at once a law of large numbers, that is, $\frac{X_t}{t}$ converges $P_0$-a.s. to a deterministic velocity, which is possibly 0. Our strategy is inspired by that of Zerner and Merkl in [37], where they proved an analogous zero–one law for random walks in two-dimensional i.i.d. environments. Note that Lemma 4.1 and the beginning of the proof of Theorem 4.2 are valid for all dimensions.

We will use the following notation: for every environment $\omega$, we consider two independent diffusions, called $X.$ and $Y.$. Stopping times with superscript 1, respectively 2, refer to $X.$, respectively $Y.$. We define for $\omega \in \Omega, x, y \in \mathbb{R}^d$ the product measure $P_{x,y}^\omega = P_{x,\omega} \times P_{y,\omega}$ as well as $P_{x,y} = \mathbb{E}P_{x,y}^\omega$. We recall that the first entrance time in a set $B$ is called $H_B$; see above (2.2). For every $\omega \in \Omega, x \in \mathbb{R}^d$ we write

$$(4.1) \qquad r(x,\omega) = P_{x,\omega}(l \cdot X_t \to \infty).$$

The basic idea is to first show that under the assumption $P_0(A_l \cup A_{-l}) = 1$, the two diffusions starting, respectively, in 0 and $y_L$, with $l \cdot y_L$ large, are unlikely to visit a same small ball located between their starting points; see Lemma 4.1. On the other hand, when $d = 2$, if we assume that $P_0(A_l)P_0(A_{-l}) > 0$, we can choose $y_L$ such that for large $L$, the two diffusions intersect "between 0 and $y_L$" with nonvanishing probability, thus leading to a contradiction; see Theorem 4.2.

Lemma 4.1 relies on the fact that for every $\omega$, $r(X_t, \omega)$ and $r(Y_t, \omega)$ are $P_{x,\omega}$-martingales by the Markov property, and they converge to $\mathbb{1}_{A_l}$, $P_{x,\omega}$-a.s. Loosely speaking, $X.$ and $Y.$ cannot meet in a region between their respective starting points if they are far apart, because $r(X_t, \omega)$ and $r(Y_t, \omega)$ would have to approach 1, respectively 0, in the same region.

LEMMA 4.1 ($d \geq 1$). *Consider $l \in S^{d-1}$ and assume $P_0(A_l \cup A_{-l}) = 1$. Then for any sequence $y_L, L \geq 4R$, satisfying $l \cdot y_L \geq 3L$, we have*

$$(4.2) \quad \lim_{L \to \infty} P_{0,y_L}[\text{there exists } z \in \mathcal{S}_{(L, l \cdot y_L - L)} : H^1_{B(z,R)} < \infty, H^2_{B(z,R)} < \infty] = 0.$$

PROOF. The considered set in (4.2) is measurable because it suffices to consider a countable, dense subset of $\mathcal{S}_{(L, l \cdot y_L - L)}$ in the union, since we use entrance times into open balls. For any integer $L \geq 4R$, its probability is bounded from above by

$$(4.3) \quad \begin{aligned} &P_0\left[\exists z \in \mathcal{S}_{(L, l \cdot y_L - L)} : H^1_{B(z,R)} < \infty, \sup_{y \in B(z,R)} r(y,\omega) < \tfrac{1}{2}\right] \\ &+ P_{y_L}\left[\exists z \in \mathcal{S}_{(L, l \cdot y_L - L)} : H^2_{B(z,R)} < \infty, \sup_{y \in B(z,R)} r(y,\omega) \geq \tfrac{1}{2}\right]. \end{aligned}$$

By Harnack's inequality (see [8], page 250), we have $\inf_{y \in B(z,R)} r(y,\omega) \geq c \sup_{y \in B(z,R)} r(y,\omega)$ and since $P_0(A_l \cup A_{-l}) = 1$, the expression in (4.3) is smaller than

$$P_0\left[\exists z \in \mathcal{S}_{(L, l \cdot y_L - L)} : H^1_{B(z,R)} < \infty, \sup_{y \in B(z,R)} r(y,\omega) < \frac{1}{2}, A_l\right]$$

$$+ P_0[T_{L-R} < \infty, A_{-l}]$$

(4.4)
$$+ P_{y_L}\left[\exists z \in \mathcal{S}_{(L, l \cdot y_L - L)} : H^2_{B(z,R)} < \infty, \inf_{y \in B(z,R)} r(y, \omega) \geq \frac{c}{2}, A_{-l}\right]$$
$$+ P_{y_L}[\tilde{T}^{\text{rel}}_{-L+R} < \infty, A_l].$$

The second and last terms converge to 0 as $L \to \infty$. The first term is smaller than

(4.5)
$$P_0\left[H^1_{\{z \in \mathcal{S}_{(L-R, l \cdot y_L - L + R)} : r(z, \omega) < 1/2\}} < \infty, A_l, T_{L-R} > \frac{L-R}{2\bar{K}}\right]$$
$$+ P_0\left[T_{L-R} \leq \frac{L-R}{2\bar{K}}\right],$$

where $\bar{K}$, defined in (1.2), denotes a uniform bound on the drift. From the martingale convergence theorem we know that $\lim_{t \to \infty} r(X_t, \omega) = \mathbb{1}_{A_l}$, $P_{0,\omega}$-a.s. This implies that the first term of (4.5) tends to 0 as $L \to \infty$ since $P_{0,\omega}$-a.s.,

$$\left\{\frac{L-R}{2\bar{K}} < H^1_{\{z \in \mathcal{S}_{(L-R, l \cdot y_L - L + R)} : r(z, \omega) < 1/2\}} < \infty, A_l\right\}$$
$$\subset \left\{\inf_{s \geq (L-R/2\bar{K})} r(X_s, \omega) < \frac{1}{2}, \lim_{s \to \infty} r(X_s, \omega) = 1\right\}.$$

The second term of (4.5) tends to 0 by Bernstein's inequality (see [2], Proposition 8.1, page 23). Using translation invariance, the third term in (4.4) is treated similarly. $\square$

THEOREM 4.2 ($d = 2$). *For any direction $l \in S^1$,*

(4.6) $$P_0[A_l] \in \{0, 1\}.$$

PROOF. Assume by contradiction that $P_0(A_l)P_0(A_{-l}) > 0$. For any integer $L \geq 4R$, we denote with $\Gamma_L$ the probability in (4.2) and recall that $R$ is defined in (1.5). We claim that there exists a sequence $y_L \geq 3L$, with $L \geq 4R$ such that

(4.7) $$\liminf_{L \to \infty} \Gamma_L > 0.$$

This with (4.2) yields a contradiction and Theorem 4.2 will follow. We already specify that $y_L \cdot l = 3L + 22R$. The component orthogonal to $l$ will be chosen in Lemma 4.3 below; see (4.28) and (4.31). In the first step we will use independence to separate the inner slab $\mathcal{IS}_L \stackrel{\text{def}}{=} \mathcal{S}_{(L+13/2R, l \cdot y_L - L - 13/2R)}$ from the half-spaces $\{x \in \mathbb{R}^d : x \cdot l \leq L\}$ and $\{x \in \mathbb{R}^d : x \cdot l \geq l \cdot y_L - L\}$. To achieve this, we use the coupling measure $\hat{P}^\omega_x$ starting at $x$ for the direction $l$ on the enlarged path-space $C(\mathbb{R}_+, \mathbb{R}^d) \times \{0, 1\}^{\mathbb{N}}$; see Section 3.1. For the direction $-l$, we denote the coupling measure starting at $y$ with $\tilde{P}_{y,\omega}$. We introduce the product measure $\hat{P}^\omega_{x,y} = \hat{P}_{x,\omega} \times \tilde{P}_{y,\omega}$. The Bernoulli variables, respectively, associated with $X_\cdot$ and $Y_\cdot$ are called $\lambda^1_\cdot$ and $\lambda^2_\cdot$, and $\hat{P}_{x,\omega}(\lambda^1 = 1) = \tilde{P}_{y,\omega}(\lambda^2 = 1) = \varepsilon$. For any $L \geq 4R$, we define the events

$$
(4.8) \quad \begin{aligned}
D^1 &= \left\{ T^{\text{rel},1}_L < \tilde{T}^{\text{rel},1}_{-R}, \sup_{T^{\text{rel},1}_L \leq s \leq \lceil T^{\text{rel},1}_L \rceil} |X_s - X_{T^{\text{rel},1}_L}| \leq \frac{R}{2} \right\}, \\
D^2 &= \left\{ \tilde{T}^{\text{rel},2}_{-L} < T^{\text{rel},2}_R, \sup_{\tilde{T}^{\text{rel},2}_{-L} \leq s \leq \lceil \tilde{T}^{\text{rel},2}_{-L} \rceil} |Y_s - Y_{\tilde{T}^{\text{rel},2}_{-L}}| \leq \frac{R}{2} \right\},
\end{aligned}
$$

and recall that $T_{\mathcal{IS}_L}$ denotes the exit time from $\mathcal{IS}_L$. For any $L \geq 4R$, we have the following lower bound for $\Gamma_L$ obtained by controlling the trajectories of $X_\cdot$ and $Y_\cdot$ in a symmetric way before we (almost surely) send them into the inner slab $\mathcal{IS}_L$ by requiring $\lambda^1_{\lceil T^1_L \rceil}, \lambda^2_{\lceil \tilde{T}^{\text{rel},2}_{-L} \rceil}$ to equal 1 [cf. (3.6)]:

$$
\Gamma_L \geq \mathbb{E}\hat{P}^\omega_{0,y_L}[D^1, \lambda^1_{\lceil T^1_L \rceil} = 1, D^2, \lambda^2_{\lceil \tilde{T}^{\text{rel},2}_{-L} \rceil} = 1, \exists z \in \mathcal{IS}_L,
$$

$$
(4.9) \quad H^1_{B(z,R)} \circ \hat{\theta}_{\lceil T^1_L \rceil + 1} < T^1_{\mathcal{IS}_L} \circ \hat{\theta}_{\lceil T^1_L \rceil + 1},
$$

$$
H^2_{B(z,R)} \circ \hat{\theta}_{\lceil \tilde{T}^{\text{rel},2}_{-L} \rceil + 1} < T^2_{\mathcal{IS}_L} \circ \hat{\theta}_{\lceil \tilde{T}^{\text{rel},2}_{-L} \rceil + 1}].
$$

With property (3.4), the latter expression equals

$$
(4.10) \quad \varepsilon^2 \mathbb{E}\hat{E}^\omega_{0,y_L}[D^1, D^2, g(\omega, X_{\lceil T^1_L \rceil}, Y_{\lceil \tilde{T}^{\text{rel},2}_{-L} \rceil})],
$$

where for $\omega \in \Omega, u, v \in \mathbb{R}^d$, we have defined

$$
g(\omega, u, v)
$$

$$
(4.11) \quad = \hat{P}^{\lambda^1_0=1}_{u,\omega} \times \tilde{P}^{\lambda^2_0=1}_{v,\omega}[\exists z \in \mathcal{IS}_L, H^1_{B(z,R)} \circ \hat{\theta}_1 < T^1_{\mathcal{IS}_L} \circ \hat{\theta}_1,
$$

$$
H^2_{B(z,R)} \circ \hat{\theta}_1 < T^2_{\mathcal{IS}_L} \circ \hat{\theta}_1].
$$

Using the fact that under $\hat{P}^1_{u,\omega}$, $X_1$ is uniformly distributed on the ball $B^u = B(u + 9Rl, R)$, and accordingly under $\tilde{P}^1_{v,\omega}$, $Y_1$ is uniformly distributed on the ball $\tilde{B}^v \stackrel{\text{def}}{=} B(v - 9Rl, R)$ [cf. (3.6)], we obtain from (4.11), for any $\omega \in \Omega, u, v \in \mathbb{R}^d$,

$$
(4.12) \quad g(\omega, u, v) = \frac{1}{|B(0,R)|^2} \iint h(\omega, x, y) \mathbb{1}_{\{x \in B^u\}} \mathbb{1}_{\{y \in \tilde{B}^v\}} \, dx \, dy,
$$

where for $\omega \in \Omega, x, y \in \mathbb{R}^d$, we have defined

$$
(4.13) \quad h(\omega, x, y) = P^\omega_{x,y}[\exists z \in \mathcal{IS}_L, H^1_{B(z,R)} < T^1_{\mathcal{IS}_L}, H^2_{B(z,R)} < T^2_{\mathcal{IS}_L}].
$$

For $x \in \bar{\mathcal{S}}_{(L+(15/2)R, L+(21/2)R)}$ and $y \in \bar{\mathcal{S}}_{(l \cdot y_L - L - (21/2)R, l \cdot y_L - L - (15/2)R)}$, $\omega \mapsto h(\omega, x, y)$ is $\mathcal{H}_{\mathcal{JS}_L} \subset \mathcal{H}_{\bar{\mathcal{S}}_{(l \cdot x - 4R, l \cdot y + 4R)}}$-measurable. On the other hand, the map

$$\omega \mapsto P_{0, y_L}^\omega[D^1, D^2, x \in B^{X_{\lceil T_L^1 \rceil}}, y \in \tilde{B}^{Y_{\lceil \tilde{T}_{-L}^{\text{rel},2} \rceil}}]$$

is $\mathcal{H}_{\{z \in \mathbb{R}^d \,:\, l \cdot z \leq l \cdot x - 7R\} \cup \{z \in \mathbb{R}^d \,:\, l \cdot z \geq l \cdot y + 7R\}}$ measurable. Hence, when we insert (4.12) into (4.10), finite range dependence [see (1.5)] yields

$$\begin{aligned}
(4.14) \quad \Gamma_L &\geq \frac{\varepsilon^2}{|B(0,R)|^2} \\
&\quad \times \iint P_0[D^1, x \in B^{X_{\lceil T_L^1 \rceil}}] P_{y_L}[D^2, y \in \tilde{B}^{Y_{\lceil \tilde{T}_{-L}^{\text{rel},2} \rceil}}] \mathbb{E} h(\omega, x, y) \, dx \, dy,
\end{aligned}$$

where the double integral in fact is only over $\bar{\mathcal{S}}_{(L+(15/2)R, L+(21/2)R)} \times \bar{\mathcal{S}}_{(l \cdot y_L - L - (21/2)R, l \cdot y_L - L - (15/2)R)}$. This stems from the definition of $B^u$ and $\tilde{B}^v$, and the fact that on the event $D^1$, $L - \frac{R}{2} \leq l \cdot X_{\lceil T_L^1 \rceil} \leq L + \frac{R}{2}$ and similarly on $D^2$ for $Y_{\lceil \tilde{T}_{-L}^{\text{rel},2} \rceil}$. Observe that for any $x \in \bar{\mathcal{S}}_{(L+(15/2)R, L+(21/2)R)}, y \in \bar{\mathcal{S}}_{(l \cdot y_L - L - (21/2)R, l \cdot y_L - L - (15/2)R)}$, we have

$$(4.15) \quad \mathbb{E} h(\omega, x, y) = 1 - \mathbb{E} P_{x,y}^\omega \left[ \inf_{0 \leq s \leq T_{\mathcal{JS}_L}^1, 0 \leq t \leq T_{\mathcal{IS}_L}^2} |X_s - Y_t| \geq 2R \right].$$

Using a discretization with cubes of side-length $\frac{R}{2\sqrt{d}}$ of the sets $X_{[0, T_{\mathcal{JS}_L}^1]}$ and $Y_{[0, T_{\mathcal{IS}_L}^2]}$ and with the help of finite range dependence, we see that (4.15) is larger than

$$(4.16) \quad \tilde{h}(x, y) \stackrel{\text{def}}{=} P_x \times P_y \left[ \inf_{0 \leq s \leq T_{\mathcal{JS}_L}^1, 0 \leq t \leq T_{\mathcal{JS}_L}^2} |X_s - Y_t| < R \right].$$

In view of (4.14) and (4.16), we have thus obtained the following lower bound for the initial probability: for any $L \geq 4R$,

$$\begin{aligned}
(4.17) \quad \Gamma_L &\geq \frac{\varepsilon^2 c P_0[\tilde{T}_{-R} = \infty] P_0[T_R = \infty]}{|B(0, R)|^2} \\
&\quad \times \iint \mu_L^+(B(x - 9lR, R)) \mu_L^-(B(y - y_L + 9lR, R)) \tilde{h}(x, y) \, dx \, dy,
\end{aligned}$$

where

$$(4.18) \quad \mu_L^+(\cdot) = P_0[X_{\lceil T_L^1 \rceil} \in \cdot | D^1] \quad \text{and} \quad \mu_L^-(\cdot) = P_0[Y_{\lceil T_{-L}^2 \rceil} \in \cdot | D^2],$$

with $D^1, D^2$ defined in (4.8) and where the positive constant $c$ is a lower bound for $P_{z,\omega}[\sup_{T_L \leq s \leq \lceil T_L \rceil} |X_s - X_{T_L}| \leq \frac{R}{2}] P_{z',\omega}[\sup_{\tilde{T}_{-L} \leq s \leq \lceil \tilde{T}_{-L} \rceil} |X_s - X_{\tilde{T}_{-L}}| \leq \frac{R}{2}]$, stemming from the support theorem (see [2], page 25).

The conclusion of the proof relies on the following lemma.

LEMMA 4.3 ($d = 2$). *If $P_0(A_l)P_0(A_{-l}) > 0$, then there exists $p \in (0,1)$, such that for any integer $L \geq 4R$, there are two measurable sets $A^+, A^- \subset \mathbb{R}^2$ and a point $y_L \in \mathbb{R}^2$, with $l \cdot y_L \geq 3L$, for which*

$$\tilde{h}(x,y) \geq p \quad \text{whenever } x \in A^+ \text{ and } y \in A^- \tag{4.19}$$

*and*

$$\int_{A^+} \mu_L^+(B(x - 9lR, R))\, dx > p, \quad \int_{A^-} \mu_L^-(B(y - y_L + 9lR, R))\, dy > p. \tag{4.20}$$

$\mu_L^+, \mu_L^-$ *are defined in* (4.18).

PROOF. Choose $e_2 \in S^1$ with $e_2 \cdot l = 0$. Let $a_k \in \mathbb{R}$, $k = 1,2,3$, be respective $\frac{k}{4}$-quantiles of the "second marginal" of $\mu_L^+$, chosen to be the smallest number such that $\mu_L^+(\mathbb{R}l + (-\infty, a_k]e_2) \geq \frac{k}{4}$. Let $b_k \in \mathbb{R}$, $k = 1,2,3$, be the corresponding quantiles for $\mu_L^-$. Define $A_k = [a_{k-1}, a_k], B_k = [b_{k-1}, b_k], k = 2, 3$. Choose $i, j \in \{2, 3\}$ such that $|A_i| = \min(|A_2|, |A_3|), |B_j| = \min(|B_2|, |B_3|)$. We define for integer $L \geq 4R$,

$$A^+ = [L + \tfrac{15}{2}R, L + \tfrac{21}{2}R]l + [a_{i-1} - U, a_i + U]e_2, \tag{4.21}$$

$$A^- = [-L - \tfrac{21}{2}R, -L - \tfrac{15}{2}R]l + [b_{j-1} - U, b_j + U]e_2 + y_L, \tag{4.22}$$

where $U = \frac{\sqrt{2}}{2}R$ (half the side-length of a square fitting into a ball of radius $R$) and where we recall that $y_L \cdot e_1 = 3L + 22R$. The component $y_L \cdot e_2$ will be chosen below (4.28). It is easy to check that

$$\int_{A^+} \mu_L^+(B(x - 9lR, R))\, dx$$
$$\geq 4U^2 \mu_L^+([L - \tfrac{3}{2}R + U, L + \tfrac{3}{2}R - U]l + [a_{i-1}, a_i]e_2)$$
$$\geq U^2,$$

where we recall that the first marginal of $\mu_L^+$ is supported by $\bar{\mathcal{S}}_{(L-R/2, L+R/2)}$. The same lower bound holds for $\int_{A^-} \mu_L^-(B(y - y_L + 9lR, R))\, dy$. This proves (4.20).

We next show (4.19). Adding the following two inequalities $a_3 - a_1 \geq 2|A_i|, b_3 - b_1 \geq 2|B_j|$ yields $(a_3 + b_3) + (-a_1 - b_1) \geq 2|A_i| + 2|B_j|$. Therefore at least one of the two following inequalities must hold:

$$a_3 + b_3 \geq |A_i| + |B_j|, \quad \text{case I}, \tag{4.23}$$

$$a_1 + b_1 \leq -(|A_i| + |B_j|), \quad \text{case II}. \tag{4.24}$$

Let us now examine case I. We derive a lower bound for $\tilde{h}(x,y)$ defined in (4.16) by producing a crossing of the trajectories of $X$ and $Y$ in a way that brings into play $D^1$ and $D^2$ [defined in (4.8)]. This allows us to use the measures $u_L^+, u_L^-$ and their quantiles to estimate the crossing probability. For $L \geq 4R$, $x \in \bar{\mathcal{S}}_{(L+(15/2)R, L+(21/2)R)}$ and $y \in \bar{\mathcal{S}}_{(l \cdot y_L - L - (21/2)R, l \cdot y_L - L - (15/2)R)}$:

(4.25)
$$\tilde{h}(x,y) \geq P_x \times P_y \bigg[ D^1, X_{\lceil T_L^{\text{rel},1} \rceil} \cdot e_2 > y \cdot e_2 + R,$$
$$\sup_{s \leq 1} |X_s - X_0 - 9lRs| \circ \theta_{\lceil T_L^{\text{rel},1} \rceil} \leq \frac{R}{2},$$
$$D^2, Y_{\lceil \tilde{T}_{-L}^{\text{rel},2} \rceil} \cdot e_2 > x \cdot e_2 + R,$$
$$\sup_{s \leq 1} |Y_s - Y_0 - (-9lRs)| \circ \theta_{\lceil \tilde{T}_{-L}^{\text{rel},2} \rceil} \leq \frac{R}{2} \bigg].$$

Indeed, on the above event (see Figure 4), the set $\mathcal{HS}_x \stackrel{\text{def}}{=} [L + \frac{13}{2}R, L + \frac{21}{2}R]l + (-\infty, e_2 \cdot x + \frac{R}{2}]e_2$ is connected to the line $\{z \in \mathbb{R}^2 : l \cdot z = l \cdot y_L - L - \frac{13}{2}R\}$ by a part of the trajectory of $X.$, that leaves the slab $\mathcal{IS}_L = \mathcal{S}_{(L+(13/2)R, y_L - L - (13/2)R)}$ through the "right" boundary without entering the set $\mathcal{HS}_y \stackrel{\text{def}}{=} [y_L - L - \frac{21}{2}R, y_L - L - \frac{13}{2}R]l + (-\infty, e_2 \cdot y + \frac{R}{2}]e_2$ containing $y$. This part of the trajectory divides the set $\mathcal{IS}_L \setminus \mathcal{HS}_x$ and gives rise to two connected, unbounded components, the lower one containing $y$. As the trajectory of $Y.$ leaves the slab $\mathcal{IS}_L$ through the "left" boundary without entering $\mathcal{HS}_x$, it has to intersect the part of the $X.$-trajectory separating the two connected components.

So we can bound $\tilde{h}(x,y)$ using the conditional measures $\mu_L^+, \mu_L^-$. Indeed with the support theorem (see [2], page 25) and translation invariance, it follows from (4.25) that

$$\tilde{h}(x,y) \geq cP_0[\tilde{T}_{-R}^{\text{rel}} = \infty]P_0[T_R^{\text{rel}} = \infty]$$

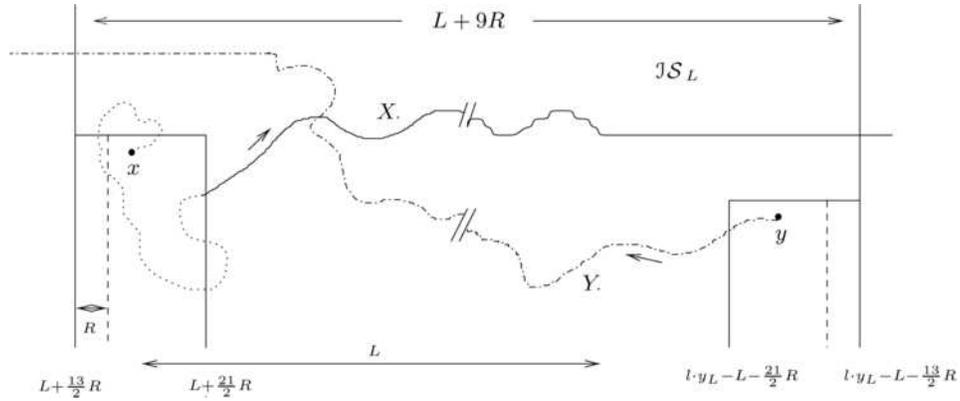

FIG. 4.

$$\times \mu_L^+(\mathbb{R}l + ((y-x)e_2 + R, \infty)e_2)\mu_L^-(\mathbb{R}l + ((x-y)e_2 + R, \infty)e_2).$$

If we choose $y_L \cdot e_2$ such that for all $x \in A^+, y \in A^-$,

(4.26) $\quad (y-x) \cdot e_2 \leq a_3 + 2U \stackrel{\text{def}}{=} \tilde{a}_3 \quad \text{and} \quad (x-y) \cdot e_2 \leq b_3 + 2U \stackrel{\text{def}}{=} \tilde{b}_3,$

then we obtain for all $x \in A^+, y \in A^-$:

$$\tilde{h}(x,y) \geq \rho \mu_L^+(\mathbb{R}l + (\tilde{a}_3 + R, \infty)e_2)\mu_L^-(\mathbb{R}l + (\tilde{b}_3 + R, \infty)e_2),$$

with $\rho = cP_0[\tilde{T}_{-R}^{\text{rel}} = \infty]P_0[T_R^{\text{rel}} = \infty] > 0$, by Lemma 3.2.

It remains to be checked that (4.26) is possible for suitable $y_L \cdot e_2$ and that

(4.27) $\quad$ (i) $\liminf_{L \to \infty} \mu_L^+(\mathbb{R}l + (\tilde{a}_3 + R, \infty)e_2) > 0,$
$\quad\quad\quad$ (ii) $\liminf_{L \to \infty} \mu_L^-(\mathbb{R}l + (\tilde{b}_3 + R, \infty)e_2) > 0.$

We first see from (4.21) and (4.22) that (4.26) is satisfied for all $x \in A^+, y \in A^-$ if

(4.28) $\quad y_L \cdot e_2 + b_j + U - a_{i-1} + U \leq \tilde{a}_3 \quad \text{and}$
$\quad\quad\quad a_i + U - y_L \cdot e_2 - b_{j-1} + U \leq \tilde{b}_3.$

Hence we have to choose $y_L \cdot e_2$ in $[a_i - b_{j-1} - b_3, a_{i-1} - b_j + a_3]$, which is possible since $(a_i - a_{i-1}) + (b_j - b_{j-1}) \leq a_3 + b_3$ in case I; see (4.23).

Finally let us check (4.27). For any $L \geq 4R$, we have [cf. (4.18)]

(4.29) $\quad \mu_L^+(\mathbb{R}l + (\tilde{a}_3 + R, \infty)e_2) \geq cP_0\left[X_{T_L} \cdot e_2 \geq \tilde{a}_3 + \frac{3R}{2}, T_L < \tilde{T}_{-R}\right],$

using the support theorem and the strong Markov property. The function $x \mapsto P_{x,\omega}[X_{T_L} \cdot e_2 \geq \tilde{a}_3 + \frac{3R}{2}, T_L < \tilde{T}_{-R}]$ is $\mathcal{L}_\omega$-harmonic in the box $(-\frac{3R}{4}, \frac{3R}{4})l + (-R, 2U + 3R)e_2$. Thus, Harnack's inequality (see [8], page 250) implies that for some constant $c > 0$, the left-hand side of (4.29) is bigger than

(4.30) $\quad \begin{aligned} &c\mathbb{E}P_{(2U+2R)e_2,\omega}\left[X_{T_L} \cdot e_2 \geq \tilde{a}_3 + \frac{3R}{2}, T_L < \tilde{T}_{-R}\right] \\ &\stackrel{\text{transl. inv.}}{=} cP_0\left[X_{T_L} \cdot e_2 \geq a_3 - \frac{R}{2}, T_L < \tilde{T}_{-R}\right] \\ &\geq cP_0[X_{\lceil T_L \rceil} \cdot e_2 \geq a_3 | D^1]P_0[D^1], \end{aligned}$

and finally the support theorem and the definition of $a_3$ yield

$$\mu_L^+(\mathbb{R}l + (\tilde{a}_3 + R, \infty)e_2) \geq cP_0[\tilde{T}_{-R} = \infty] \stackrel{\text{Lemma } 3.2}{>} 0.$$

This proves (4.27)(i). We show (4.27)(ii) in the same way. In case II [cf. (4.24)], crossings are produced by requiring instead $X_{\lceil T_L^{\text{rel},1} \rceil} \cdot e_2 < y \cdot e_2 - R$ and $Y_{\lceil \tilde{T}^{\text{rel},2}_{-L} \rceil} \cdot e_2 < x \cdot e_2 - R$ in (4.25). Moreover $y_L \cdot e_2$ has to be chosen in such a way that for all $x \in A^+, y \in A^-$,

(4.31) $\quad (y-x) \cdot e_2 \geq a_1 - 2U \quad \text{and} \quad (x-y) \cdot e_2 \geq b_1 - 2U.$

These conditions are satisfied when $y_L \cdot e_2 \in [a_1 + a_i - b_{j-1}, a_{i-1} - b_j - b_1]$, which is nonempty under (4.24). The rest of the argument has to be adjusted accordingly. This completes the proof of (4.19).

We have now proved (4.7) and as noted before Theorem 4.2 follows. $\square$

## APPENDIX

We now give the proof of Lemma 3.3.

PROOF OF LEMMA 3.3. Define the event $\Delta_0 = \{\sup_{0 \leq s \leq 1} |l \cdot (X_s - X_0)| > \frac{R}{2}\}$. The support theorem ([2], page 25) shows that there is a constant $c > 0$, such that $P_{x,\omega}(\Delta_0) < 1 - c$, for all $x \in \mathbb{R}^d$, $\omega \in \Omega$.

On the event $B_l$ [cf. (3.15)], for any $a > 0$, all the stopping times $V_k(a), k \geq 1$, are finite [recall (3.8)]. For simplicity, we drop $a$ from the notation. Define $\Delta_k = \{V_k < \infty, \sup_{V_k \leq s \leq \lceil V_k \rceil} |l \cdot (X_s - X_{V_k})| > \frac{R}{2}\}$. On the event $B_l$, $\tilde{N}_1$ is finite $P_0$-a.s., because for $n$ tending to infinity,

$$P_0\left[\bigcap_{k=1}^n \Delta_k\right] \leq \mathbb{E} E_{0,\omega}\left[\prod_{k=1}^{n-1} \mathbb{1}_{\Delta_k} P_{X_{V_n},\omega}[\Delta_0]\right] \leq (1-c)^n \stackrel{n \to \infty}{\longrightarrow} 0.$$

With the help of the strong Markov property, we obtain iteratively

$$P_0[\tilde{N}_k < \infty \text{ for all } k \geq 1 | B_l] = 1.$$

The next step is to observe that on the event $B_l$, $N_1$ is finite $\hat{P}_0$-a.s. Indeed, for any $n \geq 1$ using independence of $\lambda_j$ and $\mathcal{F}_j \otimes \mathcal{S}_{j-1}$ with respect to $\hat{P}_{x,\omega}$ [cf. (3.4)], we obtain

$$\hat{P}_x[B_l \cap \{N_1 = \infty\}]$$
$$\leq \hat{P}_0[\tilde{N}_m < \infty, \lambda_{\tilde{N}_m} = 0, \text{ for all } m \leq n]$$
$$= \sum_{j \in \mathbb{N}} \hat{P}_0[\underbrace{\tilde{N}_m < \infty, \lambda_{\tilde{N}_m} = 0, \text{ for all } m \leq n-1, \tilde{N}_n = j}_{\in \mathcal{F}_j \otimes \mathcal{S}_{j-1}}, \underbrace{\lambda_j = 0}_{\in \mathcal{S}_j}]$$
$$\stackrel{\text{induction}}{\leq} (1-\varepsilon)^n \longrightarrow 0 \quad \text{as } n \to \infty.$$

Again, by the strong Markov property, we see that on the event $B_l$, if $R_k < \infty$, then $N_{k+1} = N_1(a_k) \circ \theta_{R_k} + R_k$ is finite. $a_k$ is not time-shifted in the formula for $N_{k+1}$ [recall (3.13)]. The assumption $P_0(A_l) > 0$ and Lemma 3.2 ensure that $P_0(D = \infty) > 0$. In the next step we show that since $P_0(D =$

$\infty) > 0$, the path cannot backtrack a distance $R$ after time $N_k + 1$ for every $k \geq 1$

$$\text{(A.1)} \quad \begin{aligned} \hat{P}_0[\{R_k < \infty\} \cap B_l] &\leq \hat{P}_0[N_k < \infty, D \circ \hat{\theta}_{N_k+1} < \infty] \\ &= \mathbb{E}\hat{P}_{0,\omega}[N_k < \infty, \hat{P}^1_{X_{N_k},\omega}[\hat{P}^{\lambda_1}_{X_1,\omega}[D < \infty]]]. \end{aligned}$$

The last equality follows from (3.4). From (3.6), we see that for any $x \in \mathbb{R}^d, \omega \in \Omega$,

$$\hat{P}^1_{x,\omega}[\hat{P}^{\lambda_1}_{X_1,\omega}[D < \infty]] = \frac{1}{|B(0,R)|} \int_{B^x} P_{y,\omega}[D < \infty] \, dy.$$

Inserting this expression into (A.1), we find that for $k \geq 1$

$$\text{(A.2)} \quad \begin{aligned} &\hat{P}_0[\{R_k < \infty\} \cap B_l] \\ &\leq \frac{1}{|B(0,R)|} \int \mathbb{E}[\hat{P}_{0,\omega}[N_k < \infty, y \in B^{X_{N_k}}] P_{y,\omega}[D < \infty]] \, dy. \end{aligned}$$

The random variable $\omega \mapsto \hat{P}_{0,\omega}[N_k < \infty, y \in B^{X_{N_k}}]$ is measurable with respect to $\mathcal{H}_{\{z : z \cdot l \leq y \cdot l - 4R\}}$, because of [27], equation (3) therein, and the fact that for any $m \geq 1$, there is a $U_m \in \mathcal{F}_m \otimes \mathcal{S}_{m-1}$, with $U_m \subset \{\sup_{t \leq m} l \cdot X_t \leq l \cdot y - 7R\}$, such that $\{N_k < \infty, y \in B^{X_{N_k}}\} = \bigcup_{m \geq 1} U_m \cap \{\lambda_m = 1\}$. The random variable $\omega \mapsto P_{y,\omega}[D < \infty] = 1 - P_{0,\omega}[D = \infty]$ is measurable w.r.t. $\mathcal{H}_{\{z : z \cdot l \geq y \cdot l - R\}}$. Thus we can use the finite range dependence property (1.5) and obtain

$$\begin{aligned} \hat{P}_0[\{R_k < \infty\} \cap B_l] &\leq \hat{P}_0[N_k < \infty] P_0[D < \infty] \\ &\leq \hat{P}_0[R_{k-1} < \infty] P_0[D < \infty] \\ &\stackrel{\text{induction}}{\leq} P_0[D < \infty]^k \to 0 \qquad \text{as } k \to \infty. \end{aligned}$$

We conclude that $\hat{P}_0[\text{for some } j \geq 1 : N_j < R_j = \infty | B_l] = 1$, or in other words: $\hat{P}_0[\tau_1 < \infty | B_l] = 1$. $\square$

**Acknowledgments.** I would like to thank Prof. A.-S. Sznitman for guiding me through this work with patience and constant advice. I am also grateful for many helpful discussions with Tom Schmitz.

## REFERENCES

[1] ALON, N. and SPENCER, J. H. (1992). *The Probabilistic Method*. Wiley, New York. MR1140703
[2] BASS, R. F. (1997). *Diffusions and Elliptic Operators*. Springer, New York. MR1483890
[3] BENSOUSSAN, A., LIONS, J. L. and PAPANICOLAOU, G. (1978). *Asymptotic Analysis for Periodic Structures*. North-Holland, Amsterdam. MR0503330

[4] FRIEDMAN, A. (1975). *Stochastic Differential Equations and Applications* **1**. Academic Press, New York.
[5] COMETS, F. and ZEITOUNI, O. (2004). A law of large numbers for random walks in random mixing environment. *Ann. Probab.* **32** 880–914. MR2039946
[6] COMETS, F. and ZEITOUNI, O. (2005). Gaussian fluctuations for random walks in random mixing environments. *Israel J. Math.* **148** 87–114. MR2191225
[7] DURRETT, R. (1995). *Probability*: *Theory and Examples*, 2nd ed. Duxbury Press, Belmont, CA. MR1609153
[8] GILBARG, D. and TRUDINGER, N. S. (1998). *Elliptic Partial Differential Equations of the Second Order*. Springer, Berlin.
[9] KALIKOW, S. A. (1981). Generalized random walk in a random environment. *Ann. Probab.* **9** 753–768. MR0628871
[10] KIPNIS, C. and VARADHAN, S. R. S. (1986). A central limit theorem for additive functionals of reversible Markov processes and application to simple exclusion. *Comm. Math. Phys.* **104** 1–19. MR0834478
[11] KOMOROWSKI, T. and KRUPA, G. (2002). On the existence of invariant measure for Lagrangian velocity in compressible environments. *J. Statist. Phys.* **106** 635–651. MR1884548
[12] KOMOROWSKI, T. and KRUPA, G. (2004). On stationarity of Lagrangian observations of passive tracer velocity in a compressible environment. *Ann. Appl. Probab.* **14** 1666–1697. MR2099648
[13] KOZLOV, S. M. (1985). The method of averaging and walks in inhomogeneous environments. *Russian Math. Surveys* **40** 73–145.
[14] LANDIM, C., OLLA, S. and YAU, H. T. (1998). Convection–diffusion equation with space–time ergodic random flow. *Probab. Theory Related Fields* **112** 203–220. MR1653837
[15] LAWLER, G. F. (1982). Weak convergence of a random walk in a random environment. *Comm. Math. Phys.* **87** 81–87. MR0680649
[16] MOLCHANOV, S. (1994). Lectures on random media. *Lectures on Probability Theory. Lecture Notes in Math.* **1581** 242–411. Springer, Berlin. MR1307415
[17] OELSCHLÄGER, K. (1988). Homogenization of a diffusion process in a divergence-free random field. *Ann. Probab.* **16** 1084–1126. MR0942757
[18] OLLA, S. (1994). Homogenization of diffusion processes in random fields. Publications de l'Ecole Doctorale de l'Ecole Polytechnique, Palaiseau.
[19] OSADA, H. (1983). Homogenization of diffusion processes with random stationary coefficients. *Probability Theory and Mathematical Statistics. Lecture Notes in Math.* **1021** 507–517. Springer, Berlin. MR0736016
[20] PAPANICOLAOU, G. and VARADHAN, S. R. S. (1982). Diffusion with random coefficients. In *Statistics and Probability*: *Essays in Honor of C. R. Rao* 547–552. North-Holland, Amsterdam. MR0659505
[21] RASSOUL-AGHA, F. (2003). The point of view of the particle on the law of large numbers for random walks in a mixing random environment. *Ann. Probab.* **31** 1441–1463. MR1989439
[22] RASSOUL-AGHA, F. (2005). On the zero-one law and the law of large numbers for random walk in mixing random environment. *Electron Comm. Probab.* **10** 36–44. MR2119152
[23] RASSOUL-AGHA, F. and SEPPÄLÄINEN, T. (2005). Almost sure invariance principle for random walks in a space–time random environment. *Probab. Theory Related Fields* **133** 299–314. MR2198014

[24] RESNICK, S. I. (1992). *Adventures in Stochastic Processes*. Birkhäuser, Boston. MR1181423
[25] SCHMITZ, T. (2006). Diffusions in random environment and ballistic behavior. *Ann. Inst. H. Poincaré Probab. Statist.* To appear.
[26] SHEN, L. (2003). On ballistic random diffusions in random environment. *Ann. Inst. H. Poincaré Probab. Statist.* **39** 839–876. MR1997215
[27] SHEN, L. (2004). Addendum to the article "On ballistic random diffusions in random environment." *Ann. Inst. H. Poincaré Probab. Statist.* **40** 385–386. MR2060459
[28] STROOCK, D. and VARADHAN, S. R. S. (1979). *Multidimensional Diffusion Processes*. Springer, Berlin. MR0532498
[29] SZNITMAN, A.-S. (2000). Slowdown estimates and central limit theorem for random walks in random environment. *J. Eur. Math. Soc.* **2** 93–143. MR1763302
[30] SZNITMAN, A.-S. (2001). A class of transient random walks in random environment. *Ann. Probab.* **29** 723–764. MR1849176
[31] SZNITMAN, A.-S. (2002). An effective criterion for ballistic behaviour of random walks in random environments. *Probab. Theory Related Fields* **122** 509–544. MR1902189
[32] SZNITMAN, A.-S. and ZEITOUNI, O. (2004). On the diffusive behavior of isotropic diffusions in a random environment. *C. R. Math. Acad. Sci. Paris* **339** 429–434. MR2092758
[33] SZNITMAN, A.-S. and ZERNER, M. (1999). A law of large numbers for random walks in random environment. *Ann. Probab.* **27** 1851–1869. MR1742891
[34] VARADHAN, S. R. S. (2003). Large deviations for random walks in a random environment. *Comm. Pure Appl. Math.* **56** 1222–1245. MR1989232
[35] YURINSKY, V. V. (1980). Average of an elliptic boundary problem with random coefficients. *Siberian Math. J.* **21** 470–482.
[36] ZEITOUNI, O. (2004). Random walks in random environment. *Lectures on Probability Theory and Statistics* 183–312. *Lecture Notes in Math.* **1837** Springer, Berlin. MR2071631
[37] ZERNER, M. and MERKL, F. (2001). A zero–one law for planar random walks in random environment. *Ann. Probab.* **29** 1716–1732. MR1880239
[38] ZERNER, M. (2002). A non-ballistic law of large numbers for random walks in i.i.d. random environment. *Electron Comm. Probab.* **7** 151–197. MR1937904

DEPARTMENT OF MATHEMATICS
ETH ZURICH
CH-8092 ZURICH
SWITZERLAND
E-MAIL: goergen@math.ethz.ch